\documentclass[12pt]{imsart}

\usepackage{iftex}
\ifPDFTeX
  \usepackage[utf8]{inputenc}
  \usepackage[T1]{fontenc}
  \usepackage{lmodern}
\fi

\usepackage[margin=.75in]{geometry}
\usepackage{amsmath,amssymb,amsthm,mathtools}
\usepackage{amsfonts,mathrsfs,bbm}
\usepackage{enumitem}
\usepackage{url}
\usepackage{hyperref}
\hypersetup{colorlinks=true,citecolor=blue,linktocpage=true}
\usepackage{cleveref}
\usepackage{comment}
\usepackage{tikz}
\usetikzlibrary{arrows.meta}
\definecolor{geodesicblue}{RGB}{30,61,220}
\definecolor{geodesiccyan}{RGB}{0,148,166}
\definecolor{pointred}{RGB}{210,35,46}
\usetikzlibrary{arrows.meta,calc}
\usepackage{graphicx}
\usetikzlibrary{arrows.meta}
\RequirePackage[numbers]{natbib}
\usetikzlibrary{decorations.pathreplacing}
\startlocaldefs
\setlist[enumerate,1]{label=(\roman*)}
\setlist[enumerate,2]{label=(\alph*)}

\def\E{\mathbb{E}}
\def\P{\mathbb{P}}
\def\N{\mathbb{N}}

\def\R{\mathbb{R}}
\def\T{\mathbb T}
\def\Fbar{\bar{F}}

\newcommand{\Rinf}{\mathcal R_\infty}
\newcommand{\cS}{\mathcal S}
\newcommand{\cF}{\mathcal F}
\newcommand{\1}{\mathbf 1}
\newcommand{\ainv}{a^{\leftarrow}}

\newtheorem{dummy}{***}[section]
\newtheorem{theorem}[dummy]{Theorem}

\newtheorem{corollary}[dummy]{Corollary}
\newtheorem{definition}[dummy]{Definition}
\newtheorem{assumption}[dummy]{Assumption}
\newtheorem{lemma}[dummy]{Lemma}
\newtheorem{remark}[dummy]{Remark}
\newtheorem{proposition}[dummy]{Proposition}

\AddToHook{env/theorem/begin}{\crefalias{dummy}{theorem}}
\AddToHook{env/example/begin}{\crefalias{dummy}{example}}
\AddToHook{env/corollary/begin}{\crefalias{dummy}{corollary}}
\AddToHook{env/definition/begin}{\crefalias{dummy}{definition}}
\AddToHook{env/assumption/begin}{\crefalias{dummy}{assumption}}
\AddToHook{env/lemma/begin}{\crefalias{dummy}{lemma}}
\AddToHook{env/remark/begin}{\crefalias{dummy}{remark}}
\AddToHook{env/proposition/begin}{\crefalias{dummy}{proposition}}
\crefname{theorem}{Theorem}{Theorems}
\crefname{example}{Example}{Examples}
\crefname{corollary}{Corollary}{Corollaries}
\crefname{definition}{Definition}{Definitions}
\crefname{assumption}{Assumption}{Assumptions}
\crefname{lemma}{Lemma}{Lemmas}
\crefname{remark}{Remark}{Remarks}
\crefname{proposition}{Proposition}{Propositions}

\numberwithin{equation}{section}

\endlocaldefs
\begin{document}

\begin{frontmatter}
\title{\Large Extremal Behaviour and Geometry of Surviving Rays in Killed Branching Random Walks with Heavy-Tailed Displacements}

\begin{abstract}
For branching random walks with heavy-tailed
displacements, extremal behaviour is typically determined by
a single large jump, with little dependence on its timing.
We consider a branching random walk on a supercritical
Galton--Watson tree with killing below the superlinear
barrier $b_n=a(n/m)$, where $a$ is the inverse-tail scale
of the displacement law. While earlier work identified the
survival regime, here we investigate how the timings and
sizes of successive large jumps determine the extremal
behaviour of the alive population and surviving rays.

The cost of producing successive large jumps gives rise to two critical quantities $q<Q$, the roots of $r^2-mr+m=0$. A uniform restriction on successive new
maxima times gives upper bounds on extremal growth.
For the matching lower bounds, we introduce barrier horizons
and deficits, and develop a special-chain argument to show
that every surviving ray contains infinitely many new maxima
edges with bounded deficits. These provide suitable starting
configurations for constructing accelerating genealogies. 

Almost surely on survival, for $m>4$, the cloud and burst
speeds on the inverse-tail scale coincide and equal
$$
\gamma(m)=\frac{1+\sqrt{1-4/m}}{2}=\frac{Q}{m},
$$
whereas the sustainable speed equals $1/4$. As a corollary, we obtain the corresponding results for the branching random walk without killing.
We further establish almost sure convergence of the normalized
alive maximum, determine the exponential growth rate
$\gamma(m)$ of the alive population, and obtain a Gumbel limit for the
alive maximum with a random centering. Our assumptions cover
regularly varying and a class of log-Weibull right tails.
\end{abstract}



\runtitle{Killed Branching Random Walks}
\runauthor{Samanta \& Roy}

\author[UCL]{
\fnms{Ramkrishna Jyoti} \snm{Samanta}
\ead[label=e1]{ramkrishna.samanta.24@ucl.ac.uk}}
\author[UoL]{
\fnms{Rishideep} \snm{Roy}
\ead[label=e2]{rishideep.roy@le.ac.uk}}

\address[UCL]{
Department of Statistical Science,
University College London,
Gower Street, London WC1E 6BT, UK\printead[presep={,\ }]{e1}}
\address[UoL]{
School of Computing and Mathematical Sciences,
University of Leicester,
University Road, Leicester LE1 7RH, UK\printead[presep={,\ }]{e2}}

\end{frontmatter}
\setcounter{tocdepth}{2}
\tableofcontents
\section{Introduction}
\label{sec:introduction}
We consider a branching random walk on a supercritical Galton--Watson tree with mean offspring $\mu>1$ and heavy-tailed displacements, in which a particle is killed, together with all its descendants, as soon as its position falls below a barrier. In the absence of a killing mechanism, the exponential growth of the genealogy
creates a possibility of very large displacements. The barrier is $b_n=a(n/m)$, where $m>0$ is a parameter and $a$ is the inverse-tail scale of the displacement law, determined up to subexponential factors by $\mu^{t}\,\mathbb{P}(X\ge a(t))= \mu^{o(t)}$, where $X$ represents a generic displacement random variable. Thus $a(n)$ is the size of the largest displacement typically produced among the $\mu^n$ particles of generation $n$, and the barrier grows faster than linearly. In \citet{SamantC2026} it was shown that, under suitable constraints on the barrier and the displacement distribution, the killed process dies out almost surely for $m<4$ and survives with positive probability for $m>4$.  At $m=4$
the first-order assumptions do not determine survival: the same work
gives examples with opposite outcomes at criticality. Here we take $m>4$ and compare positions with the scale $a$: a position of order $a(sn)$ at generation $n$ is said to grow at rate $s$. On this scale the barrier grows at rate $1/m$, and $a(n)$ at rate $1$. On survival, we study three such rates: that of the maximal alive position (the \emph{cloud speed} $s_{\mathrm{cloud}}$), the largest rate attained along a subsequence of generations by an infinite surviving path (the \emph{burst speed} $s_{\mathrm{burst}}$), and the largest rate that such a path maintains at all large generations (the \emph{sustainable speed} $s_{\mathrm{sust}}$); precise definitions are given in~\eqref{eq:speeds} below. We prove that, almost surely on survival,
$$
  s_{\mathrm{cloud}}=s_{\mathrm{burst}}=\frac{1+\sqrt{1-4/m}}{2},
  \qquad
  s_{\mathrm{sust}}=\frac14 .
$$
 These speeds by their definition do not describe the same thing though: a particle attaining the generation-wise maximum may go extinct in some future generation, and an infinite surviving path that is high along a subsequence of generations may be coming closer to the barrier in between.

The survival theorem of \citet{SamantC2026} at $m>4$ shows the existence of an infinite surviving path with positive probability, but does not determine the growth rate exactly upon survival. 
 Besides the three speeds, we show that the rate of the maximal alive position converges, that the number of alive particles grows exponentially with the same exponent, and we identify the fluctuations of the maximal alive position.

\subsection{Model and main results}
\label{sec:main-results}
\label{sec:notation}

We now formally define the model starting with the branching random walk. 
Let $\T$ be a Galton--Watson tree with offspring numbers distributed as $\xi$
satisfying
\begin{equation}
\mu:=\E[\xi]\in(1,\infty).
\tag{GW1}\label{eq:GW1}
\end{equation}
We also make a further assumption that
\begin{equation}
\E[\xi^2]<\infty.
\tag{GW2}\label{eq:GW2}
\end{equation}
Write $\rho$ for the root. For vertices $v,w\in\T$,
we write $w\preceq v$ if $w$ is an ancestor of $v$
or $w=v$, and $w\prec v$ when $w\neq v$.
We say that $v$ is a descendant of $w$ if
$w\prec v$. Throughout, by the \textit{subtree below
a vertex $v$,} we mean the genealogical subtree
rooted at $v$, consisting of all its
descendants.

To each $v\in\T\setminus\{\rho\}$ attach an independent copy $X_v$ of $X$, independently of the tree. Also let $\xi_v$ be an independent copy of $\xi$ which denotes the number of children of $v$. Set $S(\rho)=0$ and
\begin{align*}
S(v):=\sum_{\rho\prec w\preceq v}X_w.
\end{align*}
The collection $\{ S(v) : v \in \T \}$ is the branching random walk. These are the positions of vertices (denoted by $v$) in the tree, which are the sums of the displacements along their ancestral lines, respectively.
For the ancestral line $\rho=v_0\prec v_1\prec\cdots\prec v_n=v$,
write
\begin{align*}
X_j(v):=X_{v_j},\qquad S_j(v):=S(v_j),\qquad 1\le j\le n.
\end{align*}
Such a $v$ is said to belong to generation $n$ of the branching random walk, and we write $|v|=n$. The set $\{w\in \T: w \preceq v\}$ denotes its ancestral lineage, including $v$ itself and the root $\rho$.
For $v\neq\rho$, let $e_v$ denote the edge entering $v$, and set
\begin{equation*}
|e_v|:=|v|,
\qquad
X_{e_v}:=X_v.
\end{equation*}

For $t\ge0$, let $\cF_t$ be the $\sigma$-field generated by $\xi_u$, $|u|<t$, and $X_v$, $1\le|v|\le t$.
We specify the barrier family $b_n^{(m)}=a(n/m)$ by choosing
a deterministic scale $a:(0,\infty)\to(0,\infty)$ relative to the right tail $\Fbar(x):=\P(X\ge x)$, subject to the following constraints:
\begin{itemize}
\item The function $a$ is strictly increasing, with $a(0)<0$, and
\begin{align}\label{C1}
\frac{a(t)}{t}\longrightarrow\infty,
\qquad t\to\infty. \tag{C1}
\end{align}

\item The scale is calibrated to the right tail in the sense that
\begin{align}\label{C2}
\log_\mu\Fbar(a(t))=-t+o(t),
\qquad t\to\infty. \tag{C2}
\end{align}

\item The scale further satisfies
\begin{align}\log_\mu\Fbar\left(\frac{a(t/r)}{t}\right)
=-\frac tr+o(t),
\qquad t\to\infty,
\quad\text{for every fixed }r>0.\label{C3}\tag{C3}\end{align}
\end{itemize}

We henceforth write $b_n=b_n^{(m)}$.
A vertex $v$ of generation $n$ is \emph{alive} when
$S_j(v)\ge b_j$ for all $1\le j\le n$. Set
\begin{equation*}
    \mathcal{A}_0 = \{ \rho \},~
    \mathcal A_n
       :=\{v\in\T:|v|=n,\ S_j(v)\ge b_j
                           \text{ for }1\le j\le n\},
    \quad A_0=1,~ A_n:=|\mathcal A_n|,
\end{equation*}
and, on $\{A_n>0\}$, define
\begin{equation*}
    M_n:=\max_{v\in\mathcal A_n}S(v).
\end{equation*}
A \emph{ray} is an infinite path $\rho=\zeta_0\prec\zeta_1\prec\cdots$ in $\mathbb{T}$ with $|\zeta_n|=n$; it is \emph{surviving} if every $\zeta_n$ is alive, and we write $S_n(\zeta):=S(\zeta_n)$.
Let $\Rinf$ denote the set of surviving 
rays and $\cS:=\{\Rinf\ne\varnothing\}$ the survival event. We
use the generalised inverse
\begin{equation*}
    \ainv(x):=\inf\{t>0:a(t)\ge x\}.
\end{equation*}

The first two quantities of interest are
\begin{equation*}
    \frac{\ainv(M_n)}n
    \qquad\text{and}\qquad
    \frac{\log_\mu A_n}{n}.
\end{equation*}
On $\cS$, $A_n\ge1$ for every $n$.
We are interested in the behavior of both quantities on the survival event $\mathcal{S}$. Here $a^{\leftarrow}(M_n)/n$ describes the height of the highest alive particle in the inverse-tail coordinate, and $\log_\mu A_n/n$ describes the exponential growth rate of the alive population.It turns out that both converge to the same limit $\gamma(m)$, defined in~\eqref{eq:gamma} below; see Theorems~\ref{thm:cloud-convergence} and~\ref{thm:alive-population-growth}.

On $\cS$ define
\begin{align}\label{eq:speeds}
    s_{\rm cloud}
       &:=\limsup_{n\to\infty}\frac{\ainv(M_n)}n, \nonumber \\
    s_{\rm burst}
       &:=\sup_{\zeta\in\Rinf}\limsup_{n\to\infty}
                      \frac{\ainv(S_n(\zeta))}{n},\\
    s_{\rm sust}
       &:=\sup_{\zeta\in\Rinf}\liminf_{n\to\infty}
                      \frac{\ainv(S_n(\zeta))}{n}. \nonumber
\end{align}
These are, respectively, the \emph{cloud}, \emph{burst} and
\emph{sustainable} speeds. 
 These quantities are measured in the inverse-tail coordinate $\ainv$, so they are rates rather than linear speeds; we keep the names cloud, burst and sustainable speed to follow \citet{BenjaminiPeres1994}, and as in \citet{Gantert2000}, where the speeds are likewise defined after a nonlinear normalization.
Clearly, as seen in \cite{Peres1999}, \cite{Lyons2005}, \cite{Gantert2000}
\begin{equation}\label{eq:speed-order}
    s_{\rm sust}\le s_{\rm burst}\le s_{\rm cloud}.
\end{equation}

The cloud and burst speeds were introduced by \citet{BenjaminiPeres1994}, while the sustainable speed goes back to \citet{LyonsPemantle1992}, who determined it, in the language of first-passage percolation, in terms of the branching number of the tree; see \citet[Section~13.8]{Lyons2005} and \citet{Peres1999}. For mean-zero increments with a finite positive exponential moment, the three speeds are positive exactly when the boundary of the tree has positive Minkowski, packing and Hausdorff dimension, respectively, and they can differ on suitable deterministic trees. On Galton--Watson trees, where these dimensions agree, the three speeds coincide \citep{Biggins1977}; see \citet[Theorem~13.32]{Lyons2005}. \citet{Gantert2000} showed that they still coincide, after a nonlinear normalisation, for semiexponential right tails. In contrast, Theorem~\ref{thm:main} shows that on a Galton--Watson tree with heavy-tailed displacements and killing, the first inequality in~\eqref{eq:speed-order} is strict.

For $m>4$ write
\begin{equation}\label{eq:gamma}
    q_\pm=q_\pm(m):=\frac{m\pm\sqrt{m^2-4m}}2,
    \qquad q:=q_-,\quad Q:=q_+,\quad
    \gamma(m):=\frac Qm=\frac1q
              =\frac{1+\sqrt{1-4/m}}2.
\end{equation}

The equalities in~\eqref{eq:gamma} use $qQ=q+Q=m$.
The proof techniques in \citet{SamantC2026} rely heavily on the analysis of the timings and magnitudes of the large displacements along the surviving rays. When one studies the growth rate along the rays, choice of the timings and magnitudes of successive jumps plays an important role. In this process we encounter the critical exponents $q_-$ and $q_+$. The interpretation of these quantities will be developed in Section~\ref{sec:mechanism}. We state the main conclusions first.

\begin{theorem}[Extremal speeds in the survival regime]
\label{thm:main}
Assume \eqref{eq:GW1}, constraints \eqref{C1}--\eqref{C3}, and $m>4$.
Almost surely on $\cS$,
\begin{equation*}
    s_{\rm cloud}\le\gamma(m),\qquad
    s_{\rm burst}\le\gamma(m),\qquad
    s_{\rm sust}\le\frac14.
\end{equation*}
If, in addition, \eqref{eq:GW2} and
Assumptions~\ref{ass:Lstar}--\ref{ass:BJ} hold, then
\begin{equation}\label{eq:main-equality}
    \P\left(\left.
       s_{\rm cloud}=s_{\rm burst}=\gamma(m),\quad
       s_{\rm sust}=\frac14\ \right|\ \cS\right)=1.
\end{equation}
\end{theorem}

The next result shows that the limsup defining $s_{\mathrm{cloud}}$ is a limit.

\begin{theorem}[Convergence of the alive cloud]
\label{thm:cloud-convergence}
Assume \eqref{eq:GW1}, \eqref{eq:GW2},  constraints \eqref{C1}--\eqref{C3},
and Assumptions~\ref{ass:Lstar}--\ref{ass:BJ}. Let $m>4$.
Then, almost surely on $\cS$,
\begin{equation*}
    \lim_{n\to\infty}\frac{\ainv(M_n)}n=\gamma(m).
\end{equation*}
\end{theorem}

\begin{theorem}[Growth of the alive population]
\label{thm:alive-population-growth}
Assume \eqref{eq:GW1}, \eqref{eq:GW2}, constraints \eqref{C1}--\eqref{C3},
and Assumptions~\ref{ass:Lstar}--\ref{ass:BJ}. Let $m>4$.
Then, almost surely on $\cS$,
\begin{equation*}
    \lim_{n\to\infty}\frac{\log_\mu A_n}{n}=\gamma(m).
\end{equation*}
\end{theorem}

In particular, the maximal position and the number of alive
particles satisfy the same leading-order relation, that is, $A_n\Fbar(M_n)= \mu^{o(n)}$. 

\begin{corollary}[Extremal speeds without killing]
\label{cor:unkilled-speeds}
Assume \eqref{eq:GW1}, \eqref{eq:GW2}, constraints
\eqref{C1}--\eqref{C3}, and
Assumptions~\ref{ass:Lstar}--\ref{ass:BJ} for arbitrarily
large $m>4$. Then, almost surely on survival of the
Galton--Watson tree, 
\begin{equation*}
    s_{\rm cloud}=s_{\rm burst}=1,
    \qquad s_{\rm sust}=\frac14.
\end{equation*}
Moreover,
\begin{equation*}
    \lim_{n\to\infty}\frac{\ainv(M_n^{(0)})}{n}
    =
    \lim_{n\to\infty}\frac{\log_\mu Z_n}{n}
    =1.
\end{equation*}
\end{corollary}

We also determine a proper centering of the maximum of the alive cloud, and convergence of the fluctuation around that to a Gumbel distribution. Let $Z_h$ be the generation-$h$ population of an unrestricted
Galton--Watson process started from one particle and write
\begin{equation}\label{eq:extreme-effective-count}
    \pi_h:=\P(Z_h>0),\qquad
    B_n^{\rm eff}
       :=\mu\sum_{j=0}^{n-1}\pi_{n-j-1}A_j.
\end{equation}
Here $\mu\pi_{n-j-1}A_j$ is the conditional expected number of children of alive vertices of generation $j$ that have an unrestricted descendant at generation $n$, and $B_n^{\rm eff}$ sums these over $j<n$.
\begin{theorem}[Convergence of the maximal alive position]
\label{thm:general-extreme-limit}
Assume \eqref{eq:GW1}, \eqref{eq:GW2}, constraints \eqref{C1}--\eqref{C3},
and Assumptions~\ref{ass:Lstar}--\ref{ass:BJ}. Let $m>4$,
and suppose that \eqref{eq:tail-regularity} holds. Put
$\Psi(x):=-\log\Fbar(x)$. Then, under $\P(\;.\;|\cS)$,
\begin{equation*}
    \Psi(M_n)-\log B_n^{\rm eff}\xrightarrow{d}G,
    \qquad \P(G\le x)=\exp\{-e^{-x}\},\quad x\in\R.
\end{equation*}
\end{theorem}
The additional tail condition in this theorem is
discussed in Section~\ref{sec:assumptions-examples}.

\begin{remark}
For a regularly varying right tail
$\Fbar(x)\sim x^{-\alpha}L(x)$ with $\alpha>0$, $L(x)$ a slowly varying function and
$\Fbar(a(t))=\mu^{-t}$, Theorem~\ref{thm:general-extreme-limit}
gives
$$
\frac{M_n}{a(\log_\mu B_n^{\rm eff})}
\xrightarrow{d}Y,
\qquad
\P(Y\le x)=\exp\{-x^{-\alpha}\},
\quad x>0.
$$
If $\Fbar(x)=\exp\{-V(\log x)\}$, the convergence takes
the form
$$
V(\log M_n)-\log B_n^{\rm eff}
\xrightarrow{d}G.
$$
\end{remark}

We conclude this section by considering the corresponding
quantities for the branching random walk without killing.
Let $Z_n$ denote the population at generation $n$ and put
$$
    M_n^{(0)}:=\max_{v\in\T:\,|v|=n}S(v).
$$
The speeds are defined as in \eqref{eq:speeds}, with the
suprema taken over all infinite rays of $\T$ and with
$M_n$ replaced by $M_n^{(0)}$.

\subsection{The mechanism behind the critical exponents}
\label{sec:mechanism}

We explain how the constants $q$ and $Q$ arise in our analysis. For simplicity, let us consider the case where the displacements are
non-negative and we want to prove that the burst speed. The simplest way we can analyze this is requiring a particle at generation $t$ to make a displacement of the scale of the barrier at $rt.$ This would mean that the expected number of  
descendants between generations $t$ and $rt$ is $\mu^{(r-1)t}$. To make sure that the burst speed is greater than $q/m$, we would require one of the descendants of the above particle near $qt$ to make a jump of the scale of the barrier at generation $q^2t$. By~\eqref{C2}, we have
$$
\P(X\ge b_{r^2t})
=
\mu^{-r^2t/m+o(t)}.
$$
Consequently, the expected number of descendants satisfying this has the following exponential order
\begin{equation*}
\mu^{\{\Phi_m(r)+o(1)\}t},
\qquad
\Phi_m(r):=
r-1-\frac{r^2}{m}
=
\frac{(r-q)(Q-r)}{m}.
\end{equation*}
Since $m>4$, we have $1<q<2<Q$, and
$$
\Phi_m(r)>0
\quad\Longleftrightarrow\quad
q<r<Q.
$$
Thus, the expected number of suitable successors grows
exponentially for $q<r<Q$, while vanishes at the roots.

We now state two results concerning the successive
new maxima times along alive particles. For an ancestral
path with increments $X_1,X_2,\ldots$, let
$$
\tau_1:=1,\qquad
\tau_{i+1}:=
\min\{j>\tau_i:X_j>X_{\tau_i}\},
$$
whenever the minimum exists. For a vertex $v$, we write $\tau_i(v)$ for these times along the ancestral line of $v$. As in \citet{SamantC2026},
if $\tau_{i+1}$ exists and the path is alive up to generation $\tau_{i+1}-1$, then
\begin{equation*}
X_{\tau_i}
\ge
\frac{b_{\tau_{i+1}-1}}{\tau_{i+1}-1}.
\end{equation*}

The following theorem gives a restriction on these times,
uniformly over all alive particles. It will be used to
establish the upper bounds in \cref{thm:main}.

\begin{theorem}[Restriction on successive new maxima times]
\label{thm:new-maxima-times}
Assume \eqref{eq:GW1}, constraints \eqref{C1}--\eqref{C3}, and
$m>4$. For every $M>m$ and $\delta>0$, almost surely
there exists a finite random time $T$ such that,
for every alive vertex $v$ and every successive new
maximum time $\tau_i(v)\ge T$, with $i\ge2$,
\begin{equation}\label{eq:new-maxima-constraint}
\frac{\tau_2(v)+\cdots+\tau_i(v)}{\tau_i(v)}
\ge q_-(M)-\delta,
\end{equation}
where
$$
q_-(M):=\frac{M-\sqrt{M^2-4M}}2.
$$
\end{theorem}

For the lower bounds, the main difficulty occurs at
the lower root $q$, where $\Phi_m(q)=0$, for the same reason as mentioned above.

For an edge $e$, define its largest barrier index by
\begin{equation}\label{eq:barrier-horizon}
H(e):=
\max\bigl(\{h\in\N:X_e\ge b_h\}\cup\{0\}\bigr),
\end{equation}
and set
\begin{equation}\label{eq:deficit}
\Delta(e):=q|e|-H(e).
\end{equation}
Here, and everywhere else we assume $\N=\{1,2,\dots\}$.
The following result shows that, along every surviving ray, infinitely many of the deficits at successive new maxima times are at most a fixed constant. 

\begin{theorem}[Bounded deficits at successive new maxima times]
\label{thm:bounded-deficit}
Assume \eqref{eq:GW1}, constraint \eqref{C1},
Assumptions~\ref{ass:u-l-bound} and \ref{ass:BJ},
and let $m>4$. There exists a deterministic $K<\infty$
such that, almost surely, simultaneously for every surviving ray, infinitely many of its edges $e$
at successive new maxima times satisfy

\begin{equation*}
\Delta(e)\le K,
\qquad\text{equivalently}\qquad
H(e)\ge q|e|-K.
\end{equation*}
\end{theorem}
The result is particularly useful for crossing the lower root, $q$ threshold.

\subsection{Assumptions and examples}
\label{sec:assumptions-examples}
\label{sec:examples}

We mention the assumptions and their roles. Constraints \eqref{C1}--\eqref{C2} specify the barrier scale: \eqref{C1}
ensures $b_n/n\to\infty$, while \eqref{C2} gives
\begin{equation*}
    \Fbar(b_n)=\mu^{-n/m+o(n)}.
\end{equation*}
These assumptions imply that $X$ has no finite positive
exponential moment. For
instance, when $\Fbar$ is eventually continuous and strictly
decreasing and $-\log\Fbar(x)=o(x)$, the inverse of
$x\mapsto-\log_\mu\Fbar(x)$ provides a scale satisfying
\eqref{C1}--\eqref{C2}. 

Note that
\begin{equation*}
    \P\left(\sum_{j=1}^nX_j\ge a(\beta n)\right)
      \le n\Fbar\!\left(\frac{a(\beta n)}n\right)
      =\mu^{-\beta n+o(n)}
      \qquad(\beta>0).
\end{equation*}
Thus \eqref{C3} ensures that the probability of the
path sum exceeding $a(\beta n)$ is, on the logarithmic
scale, bounded above by the probability of a single
displacement exceeding $a(\beta n)$.

For constructing the lower-bound, we further require the following hypotheses on the displacement law.

\begin{assumption}[Left-tail summability]
\label{ass:Lstar}
For every fixed $c>0$,
\begin{equation}\tag{$L_m^*$}\label{eq:Lstar}
    \sum_{n=2}^\infty
    n\,\P\left(X<-c\frac{b_{n-1}}n\right)<\infty.
\end{equation}
\end{assumption}

This bounds the cumulative negative displacements between the
successive large jumps selected by the survival constructions.
It strengthens the condition $(L_m)$ of \citet{SamantC2026},
which is imposed for the single coefficient $c=1$.

Note that, by iterating \eqref{eq:scale-separation} in
Lemma~\ref{lem:scale-separation}, we have that the barrier $b_n$
grows faster than any polynomial. Thus, by Markov's inequality,
\eqref{eq:Lstar} holds whenever $\E[(X^-)^\eta]<\infty$
for some $\eta>0$, where $X^-=\max\{-X,0\}$.

\begin{assumption}[Two-sided tail calibration]
\label{ass:u-l-bound}
There exist $0<c_*\le C_*<\infty$ and $t_0<\infty$ such that
\begin{equation}\tag{CAL}\label{eq:u-l-bound}
    c_*\mu^{-t}\le\Fbar(a(t))\le C_*\mu^{-t},
    \qquad t\ge t_0.
\end{equation}
\end{assumption}
It is obvious that \eqref{eq:u-l-bound} implies \eqref{C2}. Unlike \eqref{C2}, \eqref{eq:u-l-bound} supplies bounds uniform in
$t$ and will be useful in providing bounds to number of offsprings of certain type and also in the accelerating construction.

\begin{assumption}[Uniform one-big-jump condition]
\label{ass:BJ}
There is an integer $L\ge 1$ such that, for iid
$X_1,\ldots,X_n$,
\begin{equation}\tag{BJ}\label{eq:BJ}
    \sum_{n>L}\mu^n
      \P\left(\sum_{i=1}^nX_i\ge b_n,\quad
                 \max_{i\le n}X_i<b_{n-L}\right)<\infty.
\end{equation}
\end{assumption}

Assumption \eqref{eq:BJ} implies, by a first-moment estimate and
Borel--Cantelli, that almost surely every sufficiently late
particle above $b_n$ has an ancestral displacement at least
$b_{n-L}$, uniformly over the tree. This is stronger than the
exponential upper bound provided by \eqref{C3}.

We assume in addition that
\begin{equation}\tag{TR}\label{eq:tail-regularity}
    \lim_{\varepsilon\downarrow0}\limsup_{x\to\infty}
       \frac{\Fbar((1-\varepsilon)x)}
            {\Fbar((1+\varepsilon)x)}=1.
\end{equation}
This will be useful in obtaining the convergence of extremes.

The assumptions are satisfied jointly by the following two
classes, subject to a mild left-tail bound. 

\begin{lemma}[Two classes satisfying the displacement assumptions]
\label{lem:examples}
Fix $m>0$. Suppose that for some $\eta>0$ and $C_-<\infty$,
\begin{equation}\tag{LT}\label{eq:polynomial-left-tail}
    \P(X\le-x)\le C_-x^{-\eta}
    \qquad\text{for all sufficiently large }x.
\end{equation}
Suppose also that $\Fbar$ is eventually continuous and strictly
decreasing, and that one of the following holds:
\begin{enumerate}
\item There exist $\alpha>0$ and a slowly varying function $L$
such that
\begin{equation*}\tag{RV}
    \Fbar(x)=x^{-\alpha}L(x)
    \qquad\text{for all sufficiently large }x.
\end{equation*}
\item There exist $\beta\in(0,1)$ and an eventually continuously
differentiable, strictly increasing function $V$ such that
\begin{equation*}\tag{LW}
    \Fbar(x)=\exp\{-V(\log x)\}
    \qquad\text{for all sufficiently large }x,
\end{equation*}
where
\begin{equation}\label{eq:LW-smooth}
    V\in RV_\beta,\qquad
    \frac{yV'(y)}{V(y)}\longrightarrow\beta.
\end{equation}
\end{enumerate}
For all sufficiently large $t$, define $a(t)$ by
\begin{equation*}
    \Fbar(a(t))=\mu^{-t},
\end{equation*}
and modify $a$ on a bounded interval, if needed, to make it
strictly increasing, with $a(0)<0$. Then \textnormal{\eqref{C1}--\eqref{C3}} and
Assumptions~\ref{ass:Lstar}--\ref{ass:BJ} hold. In particular,
if \eqref{eq:GW2} also holds and $m>4$, the equality
statements above apply.
\end{lemma}

Condition \eqref{eq:tail-regularity} also holds for both classes,
so the distributional theorem applies as well. For a regularly
varying right tail, the scale has
$\log a(t)=t\log\mu/\alpha+o(t)$; in particular, for a
deterministic binary tree and $\Fbar(x)=x^{-\alpha}$ for
$x\ge1$, one may choose
\begin{equation*}
    a(t)=2^{t/\alpha},\qquad
    b_n=2^{n/(m\alpha)},\qquad
    \ainv(x)=\alpha\log_2x \quad(x\ge1).
\end{equation*}
For the second class, writing $W=V^{-1}$ gives
$\log a(t)=W(t\log\mu)$; because $1/\beta>1$, these
barriers grow more rapidly than exponentially.

\subsection{Comparison with literature}
\label{sec:related-work}
For a branching random walk on a supercritical Galton--Watson tree whose displacements have a finite positive exponential moment, the maximal position grows linearly: almost surely on survival, $M_n/n\to v:=\sup\{s:I(s)\le\log\mu\}$, where $I$ is the Cram\'er rate function of the displacement law \citep{Hammersley1974,Kingman1975,Biggins1976,Biggins1977}; see also \citet[Theorem~13.32]{Lyons2005}. A detailed analysis further yields that the minimum (or maximum) is equal to $\frac32\log n+O(1)$ in expectation; see \citep{AddarioBerryReed2009, HuShi2009}. Later, \citet{Aidekon2013} proved that the minimum, normalised and centered by $\frac32\log n$, converges in law to a Gumbel distribution with a random shift coming from the limit of the so called derivative martingale. This is the analogue of the results of \citet{Bramson1983} and \citet{LalleySellke1987} for branching Brownian motion. Theorem~\ref{thm:general-extreme-limit} is the analogue of this result in our setting. We also obtain a Gumbel limit for $\Psi(M_n)=-\log\Fbar(M_n)$, with random centering given by $\log B_n^{\mathrm{eff}}$, which is determined by the alive population. In the context of branching Brownian motion with absorption
at a linear barrier, \citet{YangZhu2024} used the techniques
developed for studying the extremal process in the unkilled
setting. A particularly interesting feature of their result is
that, in the survival regime, the leading-order, logarithmic correction and the decoration of the
limiting extremal process remain unchanged, with only the
random intensity changing. In contrast, in our setting the
killing mechanism affects even the leading-order term of the
growth of maximum.

Killing barriers for light-tailed branching random walks are often studied taking advantage of the exponential tilting. Work in this area goes back to the barrier problem of
\citet{BigginsLubachevskyShwartzWeiss1991}, where the authors relate survival to finding a balance between population growth and the large-deviation estimates of remaining above an absorbing boundary. Detailed analysis were further obtained in
\citet{GantertHuShi2011,AidekonJaffuel2011,Jaffuel2012,LiuZhang2019}. 

For a standard branching random walk with heavy-tailed
displacements, without the killing mechanism,  \citet{Durrett1983} studied maximal positions
under regularly varying tails. Extremal-process limits and
genealogical clustering in this regime have been developed by
\citet{BhattacharyaHazraRoy2017,BhattacharyaHazraRoy2018}.
Related results for different heavy-tail classes include 
\citet{Gantert2000, DGH2023, DGH2026,   BhattacharyaDyszewskiGantertPalmowski2025};
see also \citet{Stonner2026} for one-big-jump estimates.

In our setting, we can neither use exponential tilting and the associated large deviation estimates nor rely on the principle of one big jump alone to determine the extremal behaviour. The current work takes the next step beyond the survival theorem of \citet{SamantC2026}, whose proof restricts the times of successive new maxima of the displacements along a surviving ray.

The paper is organized as follows. We outline the proofs
in Section~\ref{sec:strategy}. Section~\ref{sec:cloud-upper}
contains the proof of Theorem~\ref{thm:new-maxima-times}
and establishes the upper bounds on the cloud and burst speeds.
In Section~\ref{sec:special-chains}, we prove
Theorem~\ref{thm:bounded-deficit}, establish the lower bound
for the burst speed and determine the sustainable speed.
Section~\ref{sec:cloud-convergence} contains the proofs of
Theorems~\ref{thm:cloud-convergence} and~\ref{thm:alive-population-growth},
and Section~\ref{sec:maximum-distribution} contains the proof of
Theorem~\ref{thm:general-extreme-limit}.
In Section~\ref{sec:discussion}, we discuss some open questions
and possible extensions of the present work.

\section{Strategy of proofs}\label{sec:strategy}

We give an outline of the proofs of the main results. Recall that $q=q_-(m)$ and $Q=q_+(m)$ satisfy
\begin{equation}
qQ=m,\qquad q+Q=m,\qquad q^2=m(q-1),\qquad Q^2=m(Q-1),\qquad \frac1q=\frac Qm.
\label{eq:strategy-identities}
\end{equation}

\medskip
\noindent\textbf{Successive new maxima times and the cloud upper bound.}
Let $1=\tau_1<\tau_2<\cdots$ be the successive new maxima times of the displacements along an alive path. Before the next new maxima time $\tau_{i+1}$, all displacements are at most $X_{\tau_i}$ and therefore
\begin{equation*}
X_{\tau_i}\ge b_{\tau_{i+1}-1}/(\tau_{i+1}-1).
\end{equation*}
Using an argument similar to that in the extinction proof of
\citet{SamantC2026}, we show that, for every $M>m$, the
new maxima times satisfying $\tau_2+\cdots+\tau_{i+1}\ge M\tau_i$
cannot occur at sufficiently late times. Further, a recursion argument forces 
$$\frac{\tau_2+\ldots+\tau_i}{\tau_i}\ge q_-(M)-\delta$$ after sufficiently many late new maxima times. Finally, if the position of an alive particle at generation $n$ exceeds $a(\beta n)$ for some $\beta>\gamma(m)$, then at least one of its ancestral displacements must be sufficiently large. Using the restriction on the successive new maxima times, a first-moment estimate  gives $s_{\mathrm{cloud}}\leq\gamma(m)$, and hence also $s_{\mathrm{burst}}\leq\gamma(m)$.

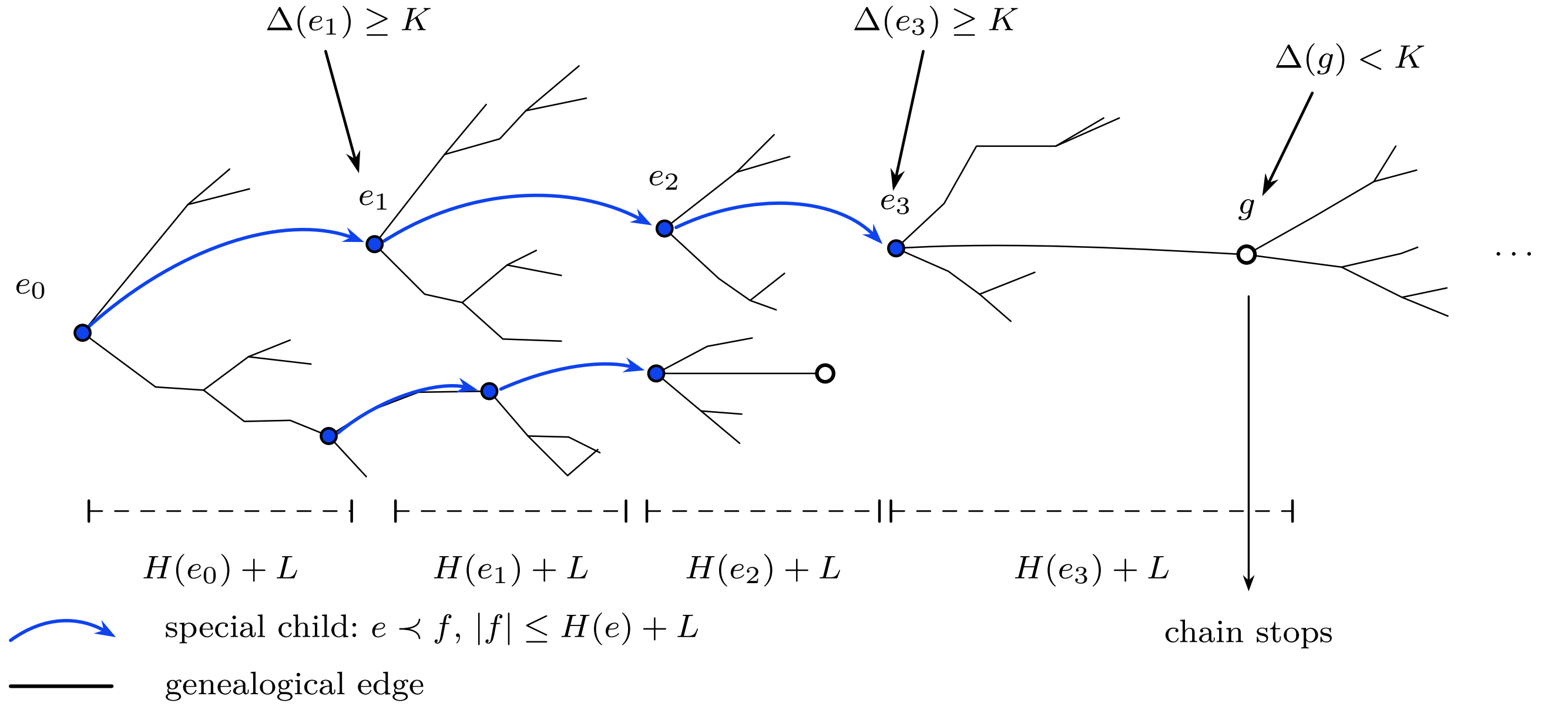
\begin{figure}[htbp]
\centering
\resizebox{\linewidth}{!}{
\begin{tikzpicture}[
x=.009cm,y=-.009cm,
line cap=round,line join=round,>=Stealth
]
\definecolor{chainblue}{RGB}{14,67,238}
\tikzset{
genealogy/.style={black,line width=.37pt},
chain/.style={chainblue,line width=.83pt,
-{Stealth[length=1.8mm,width=1.3mm]}},
horizon/.style={black,line width=.4pt,
dash pattern=on 3.3pt off 3.1pt}
}

\draw[genealogy] (200,366)--(301,243)--(341,209);
\draw[genealogy] (301,243)--(360,228);
\draw[genealogy] (200,366)--(270,418)--(316,421)
--(355,451)--(399,450)--(436,465)--(472,504);
\draw[genealogy] (316,421)--(359,389)--(399,373);
\draw[genealogy] (359,389)--(419,396);

\draw[genealogy] (480,281)--(547,195)--(587,147);
\draw[genealogy] (547,195)--(600,180)--(625,153)--(676,110);
\draw[genealogy] (625,153)--(683,141);
\draw[genealogy] (480,281)--(528,329)--(564,337)
--(603,372)--(659,374);
\draw[genealogy] (564,337)--(607,301)--(635,287);
\draw[genealogy] (607,301)--(659,311);

\draw[genealogy] (758,266)--(827,212)--(863,176);
\draw[genealogy] (827,212)--(878,197);
\draw[genealogy] (758,266)--(810,314)--(840,335)--(865,344);
\draw[genealogy] (840,335)--(873,309);

\draw[genealogy] (980,285)--(1026,242)
--(1057,187)--(1133,187)--(1179,160);
\draw[genealogy] (1133,187)--(1194,160);
\draw[genealogy] (980,285)--(1030,307)
--(1060,329)--(1113,308);
\draw[genealogy] (1060,329)--(1090,355);

\draw[genealogy] (1316,291)--(1380,255)
--(1438,221)--(1459,187);
\draw[genealogy] (1438,221)--(1479,210);
\draw[genealogy] (1316,291)--(1407,303)
--(1464,290)--(1480,284);
\draw[genealogy] (1407,303)--(1465,332)--(1509,350);
\draw[genealogy] (1465,332)--(1508,323);

\draw[genealogy] (436,465)--(475,441)--(522,423)
--(590,422)--(627,465)--(665,503);
\draw[genealogy] (627,465)--(666,466)--(696,481);
\draw[genealogy] (665,503)--(694,478);

\draw[genealogy] (750,405)--(799,379)--(842,371);
\draw[genealogy] (750,405)--(793,441)--(830,472);
\draw[genealogy] (793,441)--(832,444);
\draw[genealogy] (750,405)--(912,405);

\draw[chain] (206,360) .. controls
(297,279) and (398,252) .. (470,279);
\draw[chain] (489,278) .. controls
(577,222) and (675,225) .. (746,263);
\draw[chain] (769,265) .. controls
(842,230) and (923,236) .. (967,281);

\draw[genealogy] (980,285) .. controls
(1090,278) and (1222,286) .. (1316,291);

\draw[chain] (445,462) .. controls
(487,427) and (542,412) .. (579,421);
\draw[chain] (601,420) .. controls
(648,399) and (698,391) .. (739,402);

\foreach \xx/\yy in {
200/366,480/281,758/266,980/285,
436/465,590/422,750/405
}{
\fill[chainblue] (\xx,\yy) circle (7.4);
\draw[black,line width=.65pt]
(\xx,\yy) circle (7.4);
}

\foreach \xx/\yy in {1316/291,912/405}{
\fill[white] (\xx,\yy) circle (8);
\draw[black,line width=.9pt]
(\xx,\yy) circle (8);
}

\node[font=\scriptsize,above left=5pt]
at (200,366) {$e_0$};
\node[font=\scriptsize,above=5pt]
at (480,281) {$e_1$};
\node[font=\scriptsize,above=6pt]
at (758,266) {$e_2$};
\node[font=\scriptsize,above=5pt]
at (980,285) {$e_3$};
\node[font=\scriptsize,above=5pt]
at (1316,291) {$g$};

\node[font=\scriptsize]
at (455,68) {$\Delta(e_1)\ge K$};
\draw[black,line width=.67pt,
-{Stealth[length=1.9mm,width=1.25mm]}]
(433,96)--(465,213);

\node[font=\scriptsize]
at (1018,68) {$\Delta(e_3)\ge K$};
\draw[black,line width=.67pt,
-{Stealth[length=1.9mm,width=1.25mm]}]
(1006,96)--(977,230);

\node[font=\scriptsize]
at (1415,104) {$\Delta(g)<K$};
\draw[black,line width=.67pt,
-{Stealth[length=1.9mm,width=1.25mm]}]
(1379,136)--(1331,235);

\foreach \start/\stop/\idx in {
206/458/0,500/721/1,741/964/2,975/1360/3
}{
\draw[horizon] (\start,537)--(\stop,537);
\draw[black,line width=.65pt]
(\start,527)--(\start,547)
(\stop,527)--(\stop,547);
\node[font=\scriptsize,below=7pt]
at ({(\start+\stop)/2},537)
{$H(e_{\idx})+L$};
}

\draw[chain] (131,661) .. controls
(164,637) and (197,638) .. (232,658);
\node[font=\scriptsize,anchor=west]
at (265,650)
{special child: $e\prec f$, $|f|\le H(e)+L$};

\draw[genealogy,line width=.9pt]
(131,705)--(228,705);
\node[font=\scriptsize,anchor=west]
at (265,705) {genealogical edge};

\draw[black,line width=.5pt,
-{Stealth[length=1.4mm,width=1mm]}]
(1318,331)--(1318,614);
\node[font=\scriptsize]
at (1318,655) {chain stops};

\node[font=\scriptsize,anchor=west]
at (1539,291) {$\cdots$};

\end{tikzpicture}
}
\caption{The auxiliary $K$-special-child structure. Filled vertices are $K$-special edges. The curved arrows indicate the special-child relation, while the remaining lines represent genealogical edges. A chain stops when the next record edge has deficit less than $K$.}
\label{fig:strategy-special-chains}
\end{figure}

\medskip
\noindent\textbf{Special chains and bounded deficits.}
For a fixed $L$, assumption~\ref{ass:BJ} implies that almost surely, for every sufficiently large $n$, each particle above $b_n$ must have made displacement at least $b_{n-L}$ along its ancestral line. Take  consecutive new maxima edges, say $e_i,e_{i+1}$ along a surviving ray. This would yield
\begin{equation}
|e_{i+1}|\le H(e_i)+L.
\label{eq:strategy-horizon}
\end{equation}
Recall that $\Delta(e)=q|e|-H(e)$. We call an edge $K$-special when $\Delta(e)\ge K$, where $K\ge L$, and a descendant edge $f$ is a $K$-special child of $e$ when $f$ is $K$-special and $|e|<|f|\le H(e)+L$. The resulting descendant relation is illustrated in Figure~\ref{fig:strategy-special-chains}. If a surviving ray had only finitely many new maxima edges with deficit at most $K$, then (\ref{eq:strategy-horizon}) would give an infinite $K$-special chain along that ray.

We show that these chains are absent. To do so, we start by setting $N_j(e)$ to be the number of special children $f$ of $e$ with $j\le\Delta(f)<j+1$. For a $K$-special edge $e$ with deficit $d$, the expected number of descendants is shown to have the following bound.
\begin{equation*}
\mathbb E[N_j(e)\mid\mathcal A_v]\le C_L\mu^{-d/q+j/m},\qquad j\ge K.
\end{equation*}
For chains starting from $e_0$ whose successive deficits belong to $[j_i,j_i+1)$, we obtain the following upper bound of the form using the previous bound.
\begin{equation}
C_L^n\mu^{-\Delta(e_0)/q}\mu^{-\frac{Q-1}{m}\sum_{i=1}^{n-1}j_i+j_n/m}.
\label{eq:strategy-special-chain}
\end{equation}
Every such chain satisfies $|e_{i+1}|\le q|e_i|$, so its deficits are bounded by $C(e_0)q^i$. Next fixing $q-1<\theta<Q-1$, we show that an infinite special chain must have infinitely many indices $n$ satisfying
\begin{equation*}
\Delta(e_n)\le\theta\sum_{i=1}^{n-1}\Delta(e_i).
\end{equation*}
For such indices, the first-moment bound (\ref{eq:strategy-special-chain}) turns out to be summable when $K$ is chosen sufficiently large showing that infinite $K$-special chains can't exist. Hence, this proves that every surviving ray has infinitely many new maxima edges satisfying $\Delta(e)\le K$, which is Theorem~\ref{thm:bounded-deficit}.

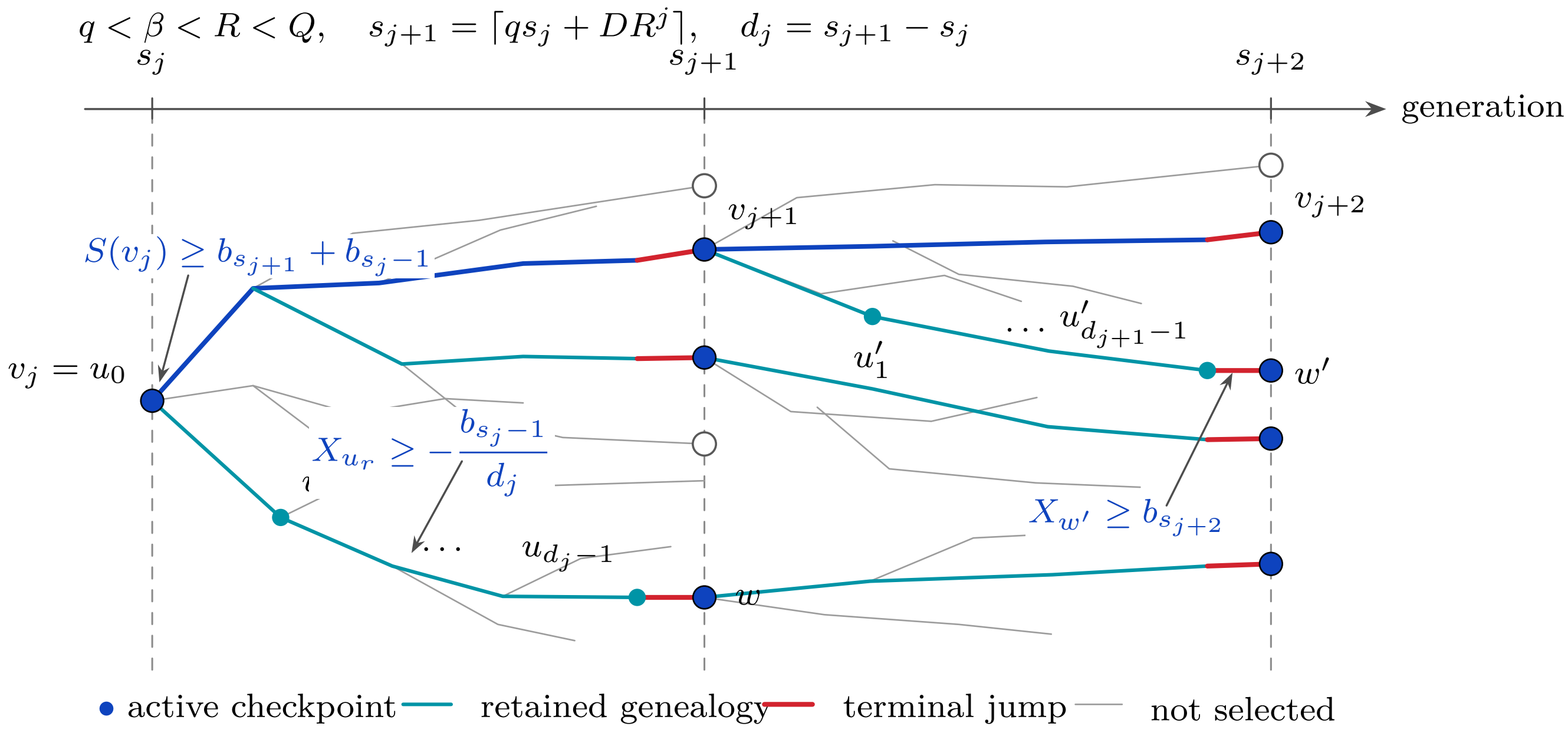
\begin{figure}[htbp]
\centering
\resizebox{\linewidth}{!}{%
\begin{tikzpicture}[
    x=1cm,y=0.9cm,
    line cap=round,line join=round,>=Stealth,
    genealogy/.style={draw=black!39,line width=.37pt},
    retained/.style={draw=geodesiccyan,line width=.85pt},
    embedded/.style={draw=geodesicblue,line width=1.1pt},
    terminal/.style={draw=pointred,line width=1.15pt},
    guide/.style={
        draw=black!45,
        dash pattern=on 2.7pt off 2.7pt,
        line width=.43pt
    },
    pointer/.style={
        draw=black!70,line width=.48pt,
        -{Stealth[length=1.4mm,width=1mm]}
    },
    eqnote/.style={
        font=\scriptsize,
        text=geodesicblue,
        fill=white,inner sep=.55pt
    },
    tinylabel/.style={font=\scriptsize,text=black}
]

\definecolor{geodesicblue}{RGB}{14,67,190}
\definecolor{geodesiccyan}{RGB}{0,148,166}
\definecolor{pointred}{RGB}{210,35,46}

\node[tinylabel,anchor=west] at (.44,6.78)
{$q<\beta<R<Q,\quad
s_{j+1}=\lceil qs_j+DR^j\rceil,\quad
d_j=s_{j+1}-s_j$};

\draw[black!70,line width=.48pt,
-{Stealth[length=1.75mm]}]
(.62,6.02)--(11.46,6.02)
node[right,tinylabel]{generation};

\foreach \x/\lab in {
    1.18/{s_j},
    5.78/{s_{j+1}},
    10.50/{s_{j+2}}
}{
    \draw[guide] (\x,.83)--(\x,5.92);
    \draw[black!70,line width=.45pt]
        (\x,5.93)--(\x,6.12);
    \node[tinylabel,anchor=south]
        at (\x,6.20) {$\lab$};
}

\draw[genealogy]
(1.18,3.32)--(2.02,3.46)--(2.85,3.20)
--(3.62,3.34)--(4.27,3.30);

\draw[genealogy]
(2.02,3.46)--(2.77,2.83)
--(3.50,2.92)--(4.20,3.15);

\draw[genealogy]
(2.02,4.36)--(2.85,4.87)
--(3.90,4.99)--(5.78,5.31);

\draw[genealogy]
(3.07,4.41)--(4.08,4.90)--(4.88,5.12);

\draw[genealogy]
(3.26,3.66)--(3.80,3.18)
--(4.60,2.98)--(5.78,2.92);

\draw[genealogy]
(2.25,2.24)--(3.15,2.73)
--(4.05,2.52)--(5.78,2.58);

\draw[genealogy]
(3.18,1.79)--(4.06,1.25)--(4.70,1.10);

\draw[genealogy]
(4.10,1.51)--(4.75,1.86)--(5.50,1.97);

\draw[genealogy]
(5.78,4.72)--(6.55,5.20)
--(7.70,5.32)--(8.80,5.30)--(10.50,5.50);

\draw[genealogy]
(7.35,4.80)--(7.90,4.49)
--(8.85,4.38)--(9.42,4.22);

\draw[genealogy]
(5.78,4.72)--(6.75,4.31)
--(7.78,4.48)--(8.42,4.24);

\draw[genealogy]
(5.78,3.72)--(6.50,3.22)
--(7.67,3.13)--(8.55,3.35);

\draw[genealogy]
(6.72,3.26)--(7.32,2.69)
--(8.54,2.56)--(9.41,2.52);

\draw[genealogy]
(5.78,1.50)--(6.78,1.34)
--(7.90,1.25)--(8.67,1.16);

\draw[genealogy]
(7.16,1.65)--(8.02,2.05)--(8.92,2.10);

\draw[embedded]
(1.18,3.32)--(2.02,4.36)
--(3.07,4.41)--(4.27,4.59)--(5.22,4.62);

\draw[retained]
(2.02,4.36)--(3.26,3.66)
--(4.27,3.73)--(5.22,3.71);

\draw[retained]
(1.18,3.32)--(2.25,2.24)
--(3.18,1.79)--(4.10,1.51)--(5.22,1.50);

\draw[terminal]
(5.22,4.62)--(5.78,4.72);

\draw[terminal]
(5.22,3.71)--(5.78,3.72);

\draw[terminal]
(5.22,1.50)--(5.78,1.50);

\draw[embedded]
(5.78,4.72)--(7.18,4.75)
--(8.62,4.79)--(9.97,4.81);

\draw[retained]
(5.78,4.72)--(7.18,4.10)
--(8.65,3.78)--(9.97,3.60);

\draw[retained]
(5.78,3.72)--(7.18,3.43)
--(8.64,3.08)--(9.97,2.96);

\draw[retained]
(5.78,1.50)--(7.16,1.65)
--(8.68,1.71)--(9.97,1.79);

\draw[terminal]
(9.97,4.81)--(10.50,4.88);

\draw[terminal]
(9.97,3.60)--(10.50,3.60);

\draw[terminal]
(9.97,2.96)--(10.50,2.97);

\draw[terminal]
(9.97,1.79)--(10.50,1.81);

\foreach \x/\y in {
    1.18/3.32,
    5.78/4.72,
    5.78/3.72,
    5.78/1.50,
    10.50/4.88,
    10.50/3.60,
    10.50/2.97,
    10.50/1.81
}{
    \fill[geodesicblue]
        (\x,\y) circle (2.7pt);
    \draw[black,line width=.37pt]
        (\x,\y) circle (2.7pt);
}

\foreach \x/\y in {
    5.78/5.31,
    5.78/2.92,
    10.50/5.50
}{
    \fill[white] (\x,\y) circle (2.75pt);
    \draw[black!65,line width=.50pt]
        (\x,\y) circle (2.75pt);
}

\fill[geodesiccyan]
    (2.25,2.24) circle (2.05pt);

\fill[geodesiccyan]
    (5.22,1.50) circle (2.05pt);

\fill[geodesiccyan]
    (7.18,4.10) circle (2.05pt);

\fill[geodesiccyan]
    (9.97,3.60) circle (2.05pt);

\node[tinylabel,anchor=east]
    at (1.10,3.55) {$v_j=u_0$};

\node[tinylabel,anchor=south west]
    at (2.30,2.30) {$u_1$};

\node[tinylabel]
    at (3.62,1.94) {$\cdots$};

\node[tinylabel,anchor=south east]
    at (5.17,1.60) {$u_{d_j-1}$};

\node[tinylabel,anchor=west]
    at (5.91,1.50) {$w$};

\node[tinylabel,anchor=south west]
    at (5.85,4.77) {$v_{j+1}$};

\node[tinylabel,anchor=south west]
    at (10.57,4.88) {$v_{j+2}$};

\node[tinylabel,anchor=north]
    at (7.18,4.01) {$u'_1$};

\node[tinylabel]
    at (8.48,3.97) {$\cdots$};

\node[tinylabel,anchor=south east]
    at (9.95,3.67)
    {$u'_{d_{j+1}-1}$};

\node[tinylabel,anchor=west]
    at (10.57,3.60) {$w'$};

\node[eqnote,anchor=west]
    at (.58,4.66)
    {$S(v_j)\ge b_{s_{j+1}}+b_{s_j-1}$};

\draw[pointer]
    (1.48,4.48)--(1.24,3.49);
\node[eqnote,anchor=west]
    at (2.48,2.84)
    {$X_{u_r}\ge
    -\dfrac{b_{s_j-1}}{d_j}$};

\draw[pointer]
    (3.76,2.76)--(3.34,1.91);

\node[eqnote,anchor=west]
    at (8.45,2.24)
    {$X_{w'}\ge b_{s_{j+2}}$};

\draw[pointer]
    (9.63,2.35)--(10.18,3.58);

\node[tinylabel,anchor=west] at (.60,.47)
{\textcolor{geodesicblue}{$\bullet$}\
 active checkpoint};

\draw[retained]
(3.27,.51)--(3.67,.51);

\node[tinylabel,anchor=west]
at (3.79,.47) {retained genealogy};

\draw[terminal]
(6.28,.51)--(6.68,.51);

\node[tinylabel,anchor=west]
at (6.81,.47) {terminal jump};

\draw[genealogy]
(8.87,.51)--(9.26,.51);

\node[tinylabel,anchor=west]
at (9.37,.47) {not selected};

\end{tikzpicture}%
}
\caption{Schematic of the accelerating genealogy. At generation $s_j$, an active vertex has position at least $b_{s_{j+1}}+b_{s_j-1}$. Intermediate displacements are bounded below, and a terminal jump of size at least $b_{s_{j+2}}$ produces an active child at generation $s_{j+1}$.}
\label{fig:strategy-accelerator}
\end{figure}

\medskip
\noindent\textbf{From bounded deficits to accelerating jumps.}
We next prove the lower bound for the burst speed. Recall the discussion from \Cref{sec:mechanism} regarding the roots of $\Phi_m(q)$. Theorem~\ref{thm:bounded-deficit} proves presence of infinite number of suitable edge configurations. Starting from such an edge, Lemma~\ref{lem:proper-start} yields that, with probability bounded away from zero, there is an alive vertex $v_0$ at a generation $s_0$ such that
\begin{equation*}
S(v_0)\ge b_{s_1}+b_{s_0-1}.
\end{equation*}
We will refer these positions as proper starts and show that these provide the proper configuration for the accelerating jumps.

Now fix $q<R<Q$ and define the checkpoint generations by
\begin{equation*}
s_{j+1}=\lceil qs_j+DR^j\rceil.
\end{equation*}
At generation $s_j$, an active particle $v_j=u_0$ has position at least $S(v_j)\ge b_{s_{j+1}}+b_{s_j-1}.$
From such a particle, we retain paths
$$
u_0,u_1,\dots,u_{d_j-1},w,
\qquad d_j=s_{j+1}-s_j,
$$
for which the intermediate displacements satisfy
$X_{u_r}\ge -b_{s_j-1}/d_j,$
and require a terminal jump at generation $s_{j+1}$ of size at least $X_w\ge b_{s_{j+2}}.$

This produces a successive active particle $v_{j+1}$ at generation $s_{j+1}$. Repeating the same construction from $v_{j+1}$ gives the next active vertex $v_{j+2}$, as illustrated in Figure~\ref{fig:strategy-accelerator}. The construction also makes sure all selected vertices remain alive.

For $D$ sufficiently large, the selected descendants form
an embedded branching process which survives with probability
bounded away from zero. This provides $s_{\mathrm{burst}}\ge Q/m=\gamma(m)$.

The sustainable speed follows from combining the extinction result of \citet{SamantC2026} and a similar analysis as above. However, the arguments and the constructed trees differ.
\medskip

\begin{figure}[t]
\centering
\resizebox{\linewidth}{!}{
\begin{tikzpicture}[x=0.82cm,y=0.8cm,>=stealth,
  line cap=round,line join=round,
  retained/.style={draw=cyan!75!black,line width=0.85pt},
  rare/.style={draw=red!85!black,line width=1.8pt},
  high/.style={circle,fill=blue!85!black,inner sep=2.15pt},
  alive/.style={circle,fill=cyan!75!black,inner sep=1.7pt}]
\draw[->,line width=0.8pt] (0.5,6.15) -- (14.9,6.15) node[right] {generation};
\foreach \x/\lab in {1.65/u,5.4/j,8.95/u_n,13.7/n}{
  \draw[densely dashed,gray!75] (\x,1.0)--(\x,6.25);
  \node[above] at (\x,6.15) {$\lab$};
}
\draw[gray!70!black,line width=1.2pt]
 (0.48,0.95) .. controls (4.2,1.05) and (6.9,1.23) .. (9.0,1.59)
 .. controls (11.0,1.93) and (12.6,2.43) .. (14.55,3.13);
\node[gray!70!black,anchor=west] at (13.55,2.42) {$b_r=a(r/m)$};
\node[high] (v) at (1.65,3.95) {};
\node[left] at (v.west) {$v$};
\node[blue!85!black,above left] at (1.63,4.25) {$S(v)\ge a(\beta u)$};
\draw[retained] (v) -- (2.6,4.1) -- (3.45,4.05) -- (4.57,3.77);
\draw[retained] (v) -- (2.6,3.75) -- (3.5,3.63) -- (4.57,3.50);
\draw[retained] (v) -- (2.6,3.45) -- (3.55,3.30) -- (4.57,3.20);
\draw[retained] (v) -- (2.6,3.17) -- (3.5,2.92) -- (4.57,2.82);
\draw[retained] (2.6,4.1) -- (3.4,4.42) -- (4.16,4.4);
\draw[retained] (3.5,3.63) -- (4.1,3.91);
\draw[retained] (3.55,3.30) -- (4.1,3.52);
\draw[gray!50,line width=0.5pt] (v) -- (2.5,2.63) -- (3.25,2.44);
\draw[gray!50,line width=0.5pt] (2.5,2.63) -- (3.20,2.75);
\foreach \x/\y in {2.6/4.1,3.45/4.05,4.57/3.77,2.6/3.75,3.5/3.63,4.57/3.5,2.6/3.45,3.55/3.3,4.57/3.2,2.6/3.17,3.5/2.92,4.57/2.82,3.4/4.42,4.16/4.4,4.1/3.91,4.1/3.52}{
 \node[alive] at (\x,\y) {};
}
\draw[rare] (4.57,3.77) -- (5.4,5.03);
\draw[rare] (4.57,3.50) -- (5.4,4.73);
\draw[rare] (4.57,3.20) -- (5.4,4.42);
\draw[rare] (4.57,2.82) -- (5.4,4.14);
\foreach \y in {5.03,4.73,4.42,4.14}{\node[high] at (5.4,\y) {};}
\node[red!80!black,anchor=east] at (4.65,5.34) {$X_e\ge a(\beta j)$};
\node[blue!85!black,anchor=west] at (5.7,5.36)
{many $\beta$-high particles};
\node at (5.45,3.72) {$\vdots$};
\draw[retained] (5.4,4.73) -- (6.35,4.42) -- (7.15,4.23) -- (8.12,3.58);
\draw[retained] (5.4,4.42) -- (6.36,4.0) -- (7.24,3.68) -- (8.10,3.16);
\draw[retained] (5.4,4.14) -- (6.34,3.61) -- (7.12,3.38) -- (8.11,2.93);
\foreach \x/\y in {6.35/4.42,7.15/4.23,8.12/3.58,6.36/4.0,7.24/3.68,8.10/3.16,6.34/3.61,7.12/3.38,8.11/2.93}{\node[alive] at (\x,\y) {};}
\draw[rare] (8.12,3.58) -- (8.95,4.94);
\draw[rare] (8.10,3.16) -- (8.95,4.44);
\draw[rare] (8.11,2.93) -- (8.95,4.13);
\foreach \y in {4.94,4.44,4.13}{\node[high] at (8.95,\y) {};}
\node[right] at (9.03,4.96) {$z$};
\node at (7.44,2.62) {$\cdots$};
\draw[retained] (8.95,4.94) -- (10.02,4.9) -- (10.93,5.10);
\draw[retained] (10.02,4.9) -- (10.98,4.56);
\draw[retained] (8.95,4.94) -- (10.03,4.27) -- (10.96,4.29);
\draw[retained] (10.03,4.27) -- (11.01,3.82);
\draw[retained] (10.93,5.10) -- (11.92,5.3) -- (12.73,5.39) -- (13.7,5.46);
\draw[retained] (11.92,5.3) -- (12.73,5.16) -- (13.7,5.04);
\draw[retained] (10.98,4.56) -- (11.93,4.72) -- (12.72,4.69) -- (13.7,4.64);
\draw[retained] (11.93,4.72) -- (12.74,4.48) -- (13.7,4.35);
\draw[retained] (10.96,4.29) -- (11.95,4.10) -- (12.74,4.00) -- (13.7,3.91);
\draw[retained] (11.95,4.10) -- (12.75,3.80) -- (13.7,3.67);
\draw[retained] (11.01,3.82) -- (11.98,3.56) -- (12.76,3.50) -- (13.7,3.40);
\draw[retained] (11.98,3.56) -- (12.77,3.34) -- (13.7,3.25);
\foreach \x/\y in {10.02/4.9,10.93/5.1,10.98/4.56,10.03/4.27,10.96/4.29,11.01/3.82,11.92/5.3,12.73/5.39,11.93/4.72,12.72/4.69,11.95/4.1,12.74/4.0,11.98/3.56,12.76/3.5}{\node[alive] at (\x,\y) {};}
\foreach \y in {5.46,5.04,4.64,4.35,3.91,3.67,3.40,3.25}{\node[alive] at (13.7,\y) {};}
\node[blue!85!black,anchor=west] at (14.02,4.57) {$A_n\ge\mu^{\vartheta(n-u_n)}$};
\draw[decorate,decoration={brace,amplitude=4pt},thick] (13.84,5.53)--(13.84,3.18);
\node[anchor=north west,scale=0.82] at (1.1,0.54) {\tikz\node[high]{};\ $\beta$-high};
\node[anchor=north west,scale=0.82] at (3.8,0.54) {\tikz\draw[retained] (0,0)--(.55,0);\ retained};
\node[anchor=north west,scale=0.82] at (6.45,0.54) {\tikz\draw[rare] (0,0)--(.55,.18);\ large displacement};
\node[anchor=north west,scale=0.82] at (10.5,0.54) {\tikz\node[alive]{};\ alive particle};
\end{tikzpicture}
}
\caption{Propagation of high particles and growth of the alive
population. Large terminal jumps produce high descendants from
which the construction is repeated. Retained descendant trees
provide the population lower bound.}
\label{fig:high-jumps-population}
\end{figure}
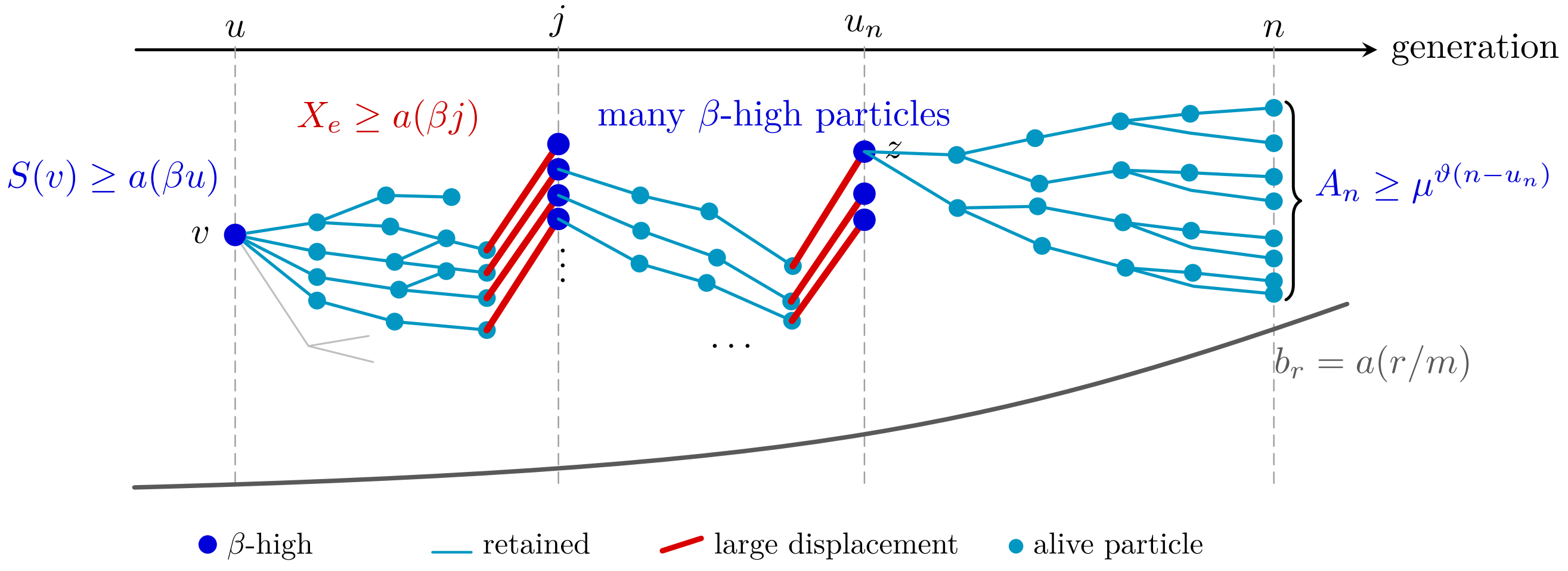

\begin{figure}[htbp]
\centering
\resizebox{\linewidth}{!}{
\begin{tikzpicture}[x=0.90cm,y=1.2cm,>=stealth,line cap=round,font=\footnotesize]
\draw[->] (0.4,1.4)--(15.5,1.4) node[right] {$n$};
\fill[red] (0.95,1.4) circle (2.2pt);
\node[left] at (0.95,1.7) {$H_t(\beta)$};
\node[below] at (0.95,1.4) {$t$};
\foreach \x in {2.65,3.1,3.65,4.25,4.85} {
\fill[red] (\x,1.4) circle (1.6pt);
\draw[blue!80!black,thin] (0.95,1.4) to[bend left=45] (\x,1.4);
}
\foreach \a/\b in {2.65/5.25,3.1/5.85,3.65/6.35,4.25/6.95} {
\fill[red] (\b,1.4) circle (1.6pt);
\draw[cyan!75!black,thin] (\a,1.4) to[bend left=43] (\b,1.4);
}
\draw (2.5,1.1)--(4.9,1.1);
\draw (2.5,1.03)--(2.5,1.17);
\draw (4.9,1.03)--(4.9,1.17);
\node[below] at (3.7,1.03) {$I_1=[At,Bt]$};
\draw (5.0,1.1)--(7.2,1.1);
\draw (5.0,1.03)--(5.0,1.17);
\draw (7.2,1.03)--(7.2,1.17);
\node[below] at (6.1,1.03) {$I_2=[L_2,U_2]$};
\draw (7.9,1.25)--(8.07,1.54);
\draw (8.2,1.25)--(8.37,1.54);
\fill[red] (9.0,1.4) circle (1.6pt);
\foreach \x in {10.2,10.85,11.5,12.2,13.05} {
\fill[red] (\x,1.4) circle (1.6pt);
\draw[blue!80!black,thin] (9.0,1.4) to[bend left=39] (\x,1.4);
}
\draw[cyan!75!black,thin] (10.2,1.4) to[bend left=38] (13.2,1.4);
\draw[cyan!75!black,thin] (10.85,1.4) to[bend left=35] (13.9,1.4);
\draw (8.8,1.1)--(12.4,1.1);
\draw (8.8,1.03)--(8.8,1.17);
\draw (12.4,1.03)--(12.4,1.17);
\node[below] at (10.6,1.03) {$I_j=[L_j,U_j]$};
\draw[red] (11.7,0.68)--(14.6,0.68);
\draw[red] (11.7,0.62)--(11.7,0.74);
\draw[red] (14.6,0.62)--(14.6,0.74);
\node[below] at (13.18,0.6) {$I_{j+1}=[L_{j+1},U_{j+1}]$};
\node[above] at (12.25,2.65) {$I_j\cap I_{j+1}\ne\varnothing$};
\end{tikzpicture}
}
\caption{Propagation of high generations. A $\beta$-high particle
produces high descendants throughout $[Au,Bu]$ on a favourable
event. Iterating the construction gives intervals which
eventually overlap and cover all sufficiently large generations.}
\label{fig:strategy-propagation}
\end{figure}
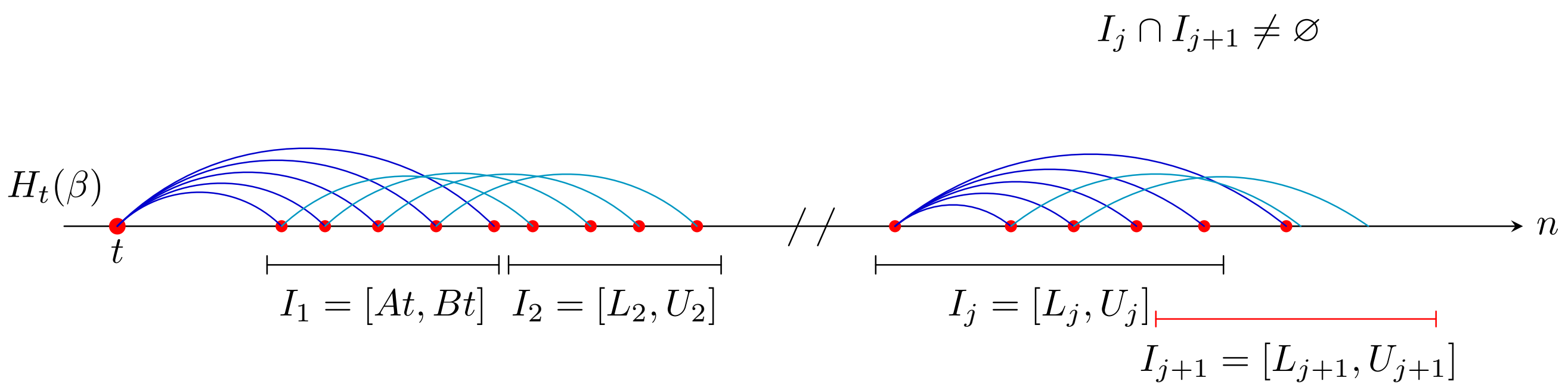

\noindent\textbf{Propagation of high particles to every generation.}
Fix $1-\gamma(m)<\beta<\gamma(m)$ and choose $\vartheta<1$
sufficiently close to one so that
$$
\frac{\vartheta}{\vartheta-\beta}<A<B<m\beta.
$$
Starting from a $\beta$-high particle $v$ at generation
$u$, i.e., $S(v)\ge a(\beta u)$, we retain descendants
with controlled negative displacements and require large
jumps $X_e\ge a(\beta j)$ at generations $j\in[Au,Bu]$.
With probability bounded away from zero, we obtain
$N_j(\beta)\ge\mu^{\delta j}$ throughout $[Au,Bu]$ for
some $\delta>0$, as illustrated in Figure~\ref{fig:high-jumps-population}.

Iterating the above construction, we obtain successive intervals
$I_1=[At,Bt]$ and $I_j=[L_j,U_j]$, which eventually
overlap and cover all sufficiently large generations
(see Figure~\ref{fig:strategy-propagation}).
Together with the burst speed, this gives
$\beta$-high particles at every sufficiently large
generation almost surely on $\cS$ gives Theorem~\ref{thm:cloud-convergence}.

For the population lower bound, we use the descendant tree of a $\beta$-high particle $z$ at
generation $u_n=\lceil n/B\rceil$ (see Figure~\ref{fig:high-jumps-population}) to show that
$A_n\ge\mu^{\vartheta(n-u_n)}$ for all sufficiently
large $n$, almost surely on $\cS$, yielding the lower bound of the population growth.

\medskip
\noindent\textbf{The distribution of the alive maximum.}
We start by considering edges whose parents
are alive and whose children have an unrestricted
descendant at generation $n$. We call such edges eligible.

Using $A_n=\mu^{\{\gamma(m)+o(1)\}n}$ on $\cS$, we cut the tree at generation
$$
s=n-\lfloor cn\rfloor,\qquad
0<c<\frac{\gamma(m)}{\gamma(m)+2}.
$$
The choice of $c$, together with the population growth
result, ensures that the probability of having
more than one displacement of the order of the maximum within any one of these descendant trees becomes negligible. This gives a Poisson approximation for
the number of eligible edges with large displacements. The analysis yields a normalization, given by
$B_n^{\mathrm{eff}},$  defined in
\eqref{eq:extreme-effective-count}.
We further show that no eligible edge before generation $s$ contribute to the maximum.

Having obtained the limiting distribution of the largest
eligible displacement $K_n$, we show that $M_n\sim K_n$ in probability under $\mathbb P(\cdot\mid\cS)$. Finally \eqref{eq:tail-regularity} yields Theorem~\ref{thm:general-extreme-limit}.

\section{Geometry of New Maxima Times}
\label{sec:cloud-upper}
Set\begin{equation*}
h_s:=\frac{b_{s-1}}{s-1},
\qquad
p_s:=\P(X\ge h_s),
\qquad s\ge2.
\end{equation*}
By \eqref{C3}, $p_s=\mu^{-s/m+o(s)}.$ In particular,
\begin{equation}
\sum_{s\ge2}p_s\mu^{\theta s}<\infty
\qquad
\text{for every }\theta<\frac1m.
\label{eq:weighted-sum}
\end{equation}
\begin{lemma}\label{lem:new-maxima-inequality}
Let $v$ be alive at generation $n$, and let
$1=\tau_1<\tau_2<\cdots<\tau_k$ be the successive new maxima times of the increments along its ancestral path. If $\tau_i<\tau_{i+1}$ are two successive such times, then
\begin{equation}
X_{\tau_i}(v)\ge h_{\tau_{i+1}}.
\label{eq:new-maxima-displacement}
\end{equation}
Moreover, every surviving ray has infinitely many successive new maxima times.
\end{lemma}

\begin{proof}
Before time $\tau_{i+1}$ every increment is at most $X_{\tau_i}(v)$, and
hence
\begin{equation*}
(\tau_{i+1}-1)X_{\tau_i}(v)
\ge
S_{\tau_{i+1}-1}(v)
\ge
b_{\tau_{i+1}-1},
\end{equation*}
which proves \eqref{eq:new-maxima-displacement}.  If a surviving ray had
only finitely many new maxima, its increments would be bounded
above by some finite constant $C$.  Then $S_n\le Cn$ for all large
$n$, contradicting \eqref{C1} and the survival condition
$S_n\ge b_n=a(n/m)$.
\end{proof}
\begin{proof}[Proof of Theorem~\ref{thm:new-maxima-times}]
For $t\ge2$, introduce the events
\begin{equation*}
E_t
:=
\left\{
\begin{array}{c}
\exists\,1=r_1<\cdots<r_k=t<r_{k+1}
\text{ and } w \text{ with } |w|=t,\\
X_{r_i}(w)\ge h_{r_{i+1}},
\qquad 1\le i\le k,\\
r_2+\cdots+r_{k+1}\ge Mt
\end{array}
\right\}.
\end{equation*}
Dropping all restrictions except the displayed displacement conditions, we obtain
\begin{equation*}
\P(E_t)
\le
\mu^t
\sum_{\substack{
A\subset\{2,3,\ldots\}\ {\rm finite}\\
\sum_{s\in A}s\ge Mt}}
\prod_{s\in A}p_s.
\end{equation*}
Choose $\eta$ such that
$1/M<\eta<1/m.$ By \eqref{eq:weighted-sum},
\begin{align*}
\P(E_t)
\le
\mu^{(1-\eta M)t}
\prod_{s\ge2}\left(1+p_s\mu^{\eta s}\right)\le
C_M\mu^{-c_Mt}
\end{align*}
for some $c_M>0$. Hence, by Borel--Cantelli, almost surely $E_t$ occurs only finitely often.

\medskip

Consider an event for which the above conclusion holds. Let
$1=\tau_1<\cdots<\tau_k$ be the successive new maxima times along an alive finite path. For
$2\le i<k$, put
\begin{equation*}
A_i:=\frac{\tau_2+\cdots+\tau_i}{\tau_i},
\qquad
q_i:=\frac{\tau_{i+1}}{\tau_i}.
\end{equation*}
For every sufficiently large new maximum time $\tau_i$, the event
$E_{\tau_i}$ does not occur. Since Lemma \ref{lem:new-maxima-inequality}
gives
$X_{\tau_j}\ge h_{\tau_{j+1}},
\; 1\le j\le i,$ we must have
\begin{equation}
A_i+q_i<M.
\label{eq:cloud-new-maxima-cost}
\end{equation}
On the other hand, $A_{i+1}
=
1+A_i/q_i.$ Therefore
\begin{align*}
A_{i+1}
&>
\frac{M}{M-A_i}
=:f(A_i).
\end{align*}

The fixed points of $f$ are $q_-(M)$ and $q_+(M)$. Moreover, we have $f(x)>x$ for
$1\le x<q_-(M).$
In particular, by continuity there exists $c=c(M,\delta)>0$ such that
\begin{equation*}
f(x)-x\ge c,
\qquad
1\le x\le q_-(M)-\delta.
\end{equation*}
Thus, whenever $A_i\le q_-(M)-\delta$ and $\tau_i$ is sufficiently
large, $A_{i+1}>A_i+c$.
Since $A_i\ge1$, there exists an integer $\ell=\ell(M,\delta)$ such
that any $\ell$ successive sufficiently late new maxima times force
$A_i\ge q_-(M)-\delta$.

It remains to show that every sufficiently long alive path contains
at least $\ell$ sufficiently late new maxima. Let $T$ be such that
\eqref{eq:cloud-new-maxima-cost} holds whenever $\tau_i\ge T$.  Suppose,
to the contrary, that there are alive vertices at arbitrarily large
generations whose ancestral lines contain fewer than $\ell$ successive
new maxima after generation $T$.
\medskip

Consider the subtree formed by all such vertices and their ancestors.
Since the Galton--Watson tree is locally finite and this subtree has
vertices at arbitrarily large generations, we get a surviving ray having at most $\ell-1$ new maxima times after generation
$T$. Hence this ray has only finitely many successive new maxima altogether,
contradicting Lemma \ref{lem:new-maxima-inequality}.

Therefore every sufficiently long alive path contains enough late
new maxima, and applying the preceding estimate to its last new maximum time $\tau_k$ gives
\begin{equation*}
\frac{\tau_2+\cdots+\tau_k}{\tau_k}
=A_k
\ge
q_-(M)-\delta.
\end{equation*}
This proves \eqref{eq:new-maxima-constraint}.
\end{proof}

\begin{proof}[Proof of the cloud upper bound in Theorem~\ref{thm:main}]
Fix $\beta>\gamma(m)$ and choose $\beta_0$ with $\gamma(m)<\beta_0<\beta$. Choose $M>m$ sufficiently close to $m$ and
then $\delta>0$ sufficiently small that, with
$c:=q_-(M)-\delta$,
one can choose $\theta<1/m$ for which
$d:=1-\theta c<\beta_0.$ Indeed, by choosing $(M,\delta,\theta)$ sufficiently close to
$(m,0,1/m)$, respectively, we can make $1-\theta c$ arbitrarily
close to $\gamma(m)$. Fix also $\beta'$ such that $d<\beta'<\beta_0.$
By Assumption {\rm\eqref{C3}}, for all sufficiently large $n$, we have
\begin{align*}
\P\left(
X\ge\frac{a(\beta_0 n)}{n}
\right)
\le
\mu^{-\beta'n}.
\end{align*}
For $1\le t\le n$, we introduce the events
\begin{equation*}
G_{t,n}
:=
\left\{
\begin{array}{c}
\exists\,1=\tau_1<\cdots<\tau_k=t
\text{ and } w \text{ with } |w|=t,\\
X_{\tau_i}(w)\ge h_{\tau_{i+1}},
\qquad 1\le i<k,\\
\tau_2+\cdots+\tau_k\ge ct,\\
X_t(w)\ge \dfrac{a(\beta_0 n)}{n}
\end{array}
\right\}.
\end{equation*}
For a fixed choice of the new maxima times, the edges appearing in the
new maximum-support inequalities are distinct from the edge entering
generation $t$. Hence, for all sufficiently large $n$,
\begin{align*}
\P(G_{t,n})
\le
\mu^{t-\beta'n}
\sum_{\substack{
A\subset\{2,\ldots,t\}\\
\sum_{s\in A}s\ge ct}}
\prod_{s\in A}p_s&\le
\mu^{t-\beta'n-\theta ct}
\prod_{s\ge2}
\left(1+p_s\mu^{\theta s}\right)\le
C_\theta\mu^{-\beta'n+dt}.
\end{align*}
Since $t\le n$ and $d<\beta'$, we have
\begin{align*}
\sum_{n\ge1}
\P\left(
\bigcup_{t=1}^nG_{t,n}
\right)
<\infty.
\end{align*}
Therefore, by Borel--Cantelli, almost surely only finitely many $n$
satisfy $\cup_{t=1}^nG_{t,n}.$

We now show that every sufficiently late alive particle satisfying $\ainv(S(v))\ge\beta n$ gives rise to one of the events $G_{t,n}$.  Let $v\in\mathcal A_n$ satisfy this inequality. Since $\beta_0n<\beta n\le\ainv(S(v))$, we have  $S(v)>a(\beta_0 n)$.
Let $t$ be the first generation at which the largest increment along
the ancestral line of $v$ is attained. Thus
$X_t(v)
=
\max_{1\le j\le n}X_j(v)$.
Since $S(v)\le nX_t(v),$
we obtain
$X_t(v) \ge S(v)/n \ge a(\beta_0 n)/n.$

Let
$1=\tau_1<\cdots<\tau_k=t$
be the successive new maxima times up to $t$. Since $t$ is the first time the
largest increment along the path to $v$ is attained, it is the last new maximum time along this path. By Theorem~\ref{thm:new-maxima-times}, for all sufficiently large $n$, $\tau_2+\cdots+\tau_k\ge ct.$ Moreover, the new maximum-support inequalities give
\begin{align*}
X_{\tau_i}(v)\ge h_{\tau_{i+1}},
\qquad 1\le i<k.
\end{align*}
Thus $G_{t,n}$ occurs. It follows that, almost surely, for every sufficiently large $n$,
there is no $v\in\mathcal A_n$ satisfying
$\ainv(S(v))\ge\beta n.$ Hence on $\cS$,
\begin{align*}
\limsup_{n\to\infty}
\frac{\ainv(M_n)}{n}
\le\gamma(m).
\end{align*}
This proves $s_{\rm cloud}\le\gamma(m)$ and since $s_{\rm burst}\leq s_{\rm cloud}$, we also have $s_{\rm burst}\le\gamma(m)$.
\end{proof}

\section{Special chains, Bounded deficits and Speeds}
\label{sec:special-chains}
\label{sec:new-maxima-horizons}

\begin{proposition}[One-big-jump statement uniform over the tree]\label{prop:uniform-big-jump}
Assume \eqref{eq:GW1} and Assumption~\ref{ass:BJ}. Almost surely, for every sufficiently large generation $n$ and every
vertex $v$ at generation $n$,
\begin{equation}
S(v)\ge b_n
\quad\Longrightarrow\quad
\max_{w\preceq v,\ w\neq\rho}X_w\ge b_{n-L},
\label{eq:uniform-big-jump}
\end{equation}
after increasing the constant $L$ in \eqref{eq:BJ} if necessary.
\end{proposition}

\begin{proof}
\eqref{eq:BJ} and Borel--Cantelli give the result.
\end{proof}

Recall that $H(e)$ is the Largest barrier index defined in \eqref{eq:barrier-horizon}.

\begin{corollary}\label{cor:new-maxima-horizons}
Almost surely, simultaneously for every surviving ray, all
sufficiently late consecutive edges at new maxima times $e_i,e_{i+1}$,
entering generations $\tau_i<\tau_{i+1}$, satisfy
\begin{equation}
\tau_{i+1}\le H(e_i)+L.
\label{eq:new-maxima-horizon}
\end{equation}
\end{corollary}

\begin{proof}
Apply Proposition~\ref{prop:uniform-big-jump} to the vertex of the ray at
generation $\tau_{i+1}-1$.  Before the next new maximum, the maximal
increment is exactly the displacement on $e_i$.  Hence
\begin{equation}
X_{e_i}\ge b_{\tau_{i+1}-1-L}.
\end{equation}
Increasing $L$ by one if necessary yields
\eqref{eq:new-maxima-horizon}.
\end{proof}

By Assumption~\ref{ass:u-l-bound}, for small $h$ there exists constants $c$ and $C$ such that
\begin{equation}
c\mu^{-h/m}
\le
\P(H(e)\ge h)
=
\Fbar(b_h)
\le
C\mu^{-h/m},
\qquad h\ge1.
\label{eq:horizon-upper-special}
\end{equation}

Recall the deficit $\Delta(e)=q|e|-H(e)$ from \eqref{eq:deficit}. Let $L$ be the deterministic constant appearing in
Corollary \ref{cor:new-maxima-horizons}. In this section $L$ is fixed.

\subsection{The descendant structure}

We use the Ulam--Harris realisation of the Galton--Watson tree.
Let
\begin{equation*}
\mathbb U:=\bigcup_{n\ge0}\mathbb N^n,
\end{equation*}
with root $\varnothing$. For a fixed deterministic label
$v\in\mathbb U$, define
\begin{equation*}
\mathscr P_v
:=
\sigma\bigl(
\xi_u:\ u\prec v;\,
X_w:\ \varnothing\prec w\preceq v
\bigr).
\end{equation*}

The event that $v$ belongs to the Galton--Watson tree is
$\mathscr P_v$-measurable. On this event, the quantities
$H(e_v)$ and $\Delta(e_v)$ are also $\mathscr P_v$-measurable.
On the other hand, $\xi_v$ and all offspring variables and
displacements strictly below $v$ are independent of
$\mathscr P_v$.

We now define the auxiliary chains used in the first-moment argument.

\begin{definition}
\label{def:special-chain}
Let $K\ge L$ be an integer. An edge $e$ is called
$K$-special if $\Delta(e)\ge K$. Given two
$K$-special edges $e$ and $f$ at generations $t$
and $u$, respectively, we say that $f$ is a
$K$-special child of $e$ if $e\prec f$ and
\begin{equation*}
t<u\le H(e)+L.
\end{equation*}
A finite or infinite sequence of edges $(e_i)$
is called a $K$-special chain if $e_{i+1}$ is
a $K$-special child of $e_i$ for every pair of
successive edges in the sequence.
\end{definition}

We first estimate the number of special children of a fixed edge. For a
$K$-special edge $e$ and an integer $j\ge K$, let $N_j(e)$ denote
the number of $K$-special children $f$ of $e$ satisfying
\begin{equation*}
j\le \Delta(f)<j+1.
\end{equation*}

\begin{lemma}[One-step estimate]
\label{lem:special-one-step}
There exists $C_L<\infty$ such that, for every deterministic
$v\in\mathbb U$, on the event that $e_v$ is $K$-special,
\begin{equation}
\E\!\left(
N_j(e_v)\mid\mathscr P_v
\right)
\le
C_L\mu^{-\Delta(e_v)/q+j/m},
\qquad j\ge K.
\label{eq:special-one-step}
\end{equation}
\end{lemma}

\begin{proof}
Fix $v$ and work on the event that $e_v$ is $K$-special. Put
\begin{equation*}
t:=|v|,
\qquad
\delta_e:=\Delta(e_v).
\end{equation*}
Both $t$ and $\delta_e$ are determined by $\mathscr P_v$. Consider a descendant edge $f$ entering generation $u>t$.
If $j\le\Delta(f)<j+1$, then
\begin{equation*}
j\le qu-H(f)<j+1,
\end{equation*}
and therefore
$H(f)>qu-j-1.$ It follows from \eqref{eq:horizon-upper-special} that
\begin{equation*}
\P\bigl(j\le \Delta(f)<j+1\bigr)
\le
C\mu^{-qu/m+j/m}.
\end{equation*}
Conditional on $\mathscr P_v$, a special child must also satisfy $u\le H(e_v)+L=qt-\delta_e+L.$ 

Consequently using \eqref{eq:strategy-identities}, we have 
\begin{align*}
\E\!\left(
N_j(e_v)\mid\mathscr P_v
\right)
&\le
C\sum_{u=t+1}^{\lfloor qt-\delta_e+L\rfloor}
\mu^{u-t}\mu^{-qu/m+j/m}
\nonumber\\
&=
C\mu^{-t+j/m}
\sum_{u=t+1}^{\lfloor qt-\delta_e+L\rfloor}
\mu^{u/q}.
\end{align*}
Hence $\E\!\left(
N_j(e_v)\mid\mathscr P_v
\right)
\le
C_L\mu^{-t+j/m}\mu^{(qt-\delta_e)/q}
=
C_L\mu^{-\delta_e/q+j/m}.$ This proves \eqref{eq:special-one-step}.
\end{proof}

\subsection{Estimates for special chains}

Fix a deterministic label $v_0\in\mathbb U$ and write
$e_0:=e_{v_0}$. For deterministic labels
\begin{equation}
v_0\prec v_1\prec\cdots\prec v_n
\label{eq:deterministic-tuple}
\end{equation}
and integers $j_1,\ldots,j_n\ge K$, let
$A(v_0,\ldots,v_n;j_1,\ldots,j_n)$
be the event that all the vertices in
\eqref{eq:deterministic-tuple} belong to the Galton--Watson tree and
\begin{equation*}
e_0\prec e_1\prec\cdots\prec e_n
\end{equation*}
is a $K$-special chain, with
$j_i\le \Delta(e_i)<j_i+1,\quad 1\le i\le n.$ Consequently, define
\begin{equation}
N_{e_0}(j_1,\ldots,j_n)
:=
\sum_{v_0\prec v_1\prec\cdots\prec v_n}
\1_{A(v_0,\ldots,v_n;j_1,\ldots,j_n)},
\label{eq:def-chain-count}
\end{equation}
counting the number of special chains
starting from $e_0$ whose successive deficits satisfy the above constraints.

\begin{lemma}[Special-chain estimate]
\label{lem:special-chain-estimate}
For every $n\ge1$, on the event that $e_0$ is $K$-special,
\begin{equation}
\E\!\left(
N_{e_0}(j_1,\ldots,j_n)
\mid\mathscr P_{v_0}
\right)
\le
C_L^n
\mu^{-\Delta(e_0)/q}
\mu^{-\frac{Q-1}{m}\sum_{i=1}^{n-1}j_i}
\mu^{j_n/m}.
\label{eq:special-chain-estimate}
\end{equation}
\end{lemma}

\begin{proof}
For $n=1$, $N_{e_0}(j_1)=N_{j_1}(e_0),$
and \eqref{eq:special-chain-estimate} follows from
Lemma \ref{lem:special-one-step}. Suppose now that $n\ge2$. For
$v_0\prec v_1\prec\cdots\prec v_{n-1}$, the event
\begin{equation}
A(v_0,\ldots,v_{n-1};j_1,\ldots,j_{n-1})
\label{eq:partial-chain-event}
\end{equation}
is $\mathscr P_{v_{n-1}}$-measurable. Conditional on the occurence of the event \eqref{eq:partial-chain-event}, the number of special child of $e_{n-1}:= e_{v_{n-1}}$ whose
next deficit lies in $[j_n,j_n+1)$ is $N_{j_n}(e_{n-1}).$ Therefore, by Lemma \ref{lem:special-one-step},
\begin{align*}
&\E\!\left(
N_{j_n}(e_{n-1})
\mid\mathscr P_{v_{n-1}}
\right)
\le
C_L\mu^{-\Delta(e_{n-1})/q+j_n/m}.
\end{align*}
On the event \eqref{eq:partial-chain-event}, $\Delta(e_{n-1})\ge j_{n-1},$ and so
\begin{equation}
\E\!\left(
N_{j_n}(e_{n-1})
\mid\mathscr P_{v_{n-1}}
\right)
\le
C_L\mu^{-j_{n-1}/q+j_n/m}.
\label{eq:chain-continuation}
\end{equation}
Therefore, summing over all possible tuples and \eqref{eq:def-chain-count}, we get
\begin{align*}
N_{e_0}(j_1,\ldots,j_n)
&=
\sum_{v_0\prec\cdots\prec v_{n-1}}
\1_{A(v_0,\ldots,v_{n-1};
j_1,\ldots,j_{n-1})}
N_{j_n}(e_{v_{n-1}}).
\end{align*}
Therefore, \eqref{eq:chain-continuation} yields
\begin{align}
&\E\!\left(
N_{e_0}(j_1,\ldots,j_n)
\mid\mathscr P_{v_0}
\right)
\le
C_L\mu^{-j_{n-1}/q+j_n/m}
\E\!\left(
N_{e_0}(j_1,\ldots,j_{n-1})
\mid\mathscr P_{v_0}
\right).
\label{eq:special-chain-recursion}
\end{align}
Iterating \eqref{eq:special-chain-recursion}, starting from the
$n=1$ estimate, yields
\begin{align*}
\E\!\left(
N_{e_0}(j_1,\ldots,j_n)
\mid\mathscr P_{v_0}
\right)&\le
C_L^n
\mu^{-\Delta(e_0)/q+j_1/m}
\prod_{i=2}^n
\mu^{-j_{i-1}/q+j_i/m}
\nonumber\\
&=
C_L^n
\mu^{-\Delta(e_0)/q}
\mu^{(\frac1m-\frac1q)
\sum_{i=1}^{n-1}j_i}
\mu^{j_n/m}.
\end{align*}
Finally, $1/q=Q/m$ proves \eqref{eq:special-chain-estimate}.
\end{proof}

\subsubsection{Absence of infinite special chains}

We next show that there exists $K$ such that no infinite $K$-special
chains exist. We start by choosing
\begin{equation*}
q-1<\theta<Q-1.
\end{equation*}
Let
\begin{equation*}
e_0\prec e_1\prec e_2\prec\cdots
\end{equation*}
be an infinite $K$-special chain, and write
$t_i:=|e_i|,
\;\Delta_i:=\Delta(e_i).$ Since $e_{i+1}$ is a special child of $e_i$, we have
\begin{equation*}
t_{i+1}
\le
H(e_i)+L
=
qt_i-\Delta_i+L.
\end{equation*}
If $K\ge L$, then $\Delta_i\ge K$, and hence $t_{i+1}\le qt_i.$ It follows that
$t_i\le q^it_0.$ Moreover, since $H(e_i)\ge0$,
\begin{equation}
\Delta_i\le qt_i\le qt_0q^i.
\label{eq:special-deficit-upper}
\end{equation}
We say that $n\ge2$ is a \emph{special index} if
\begin{equation*}
\Delta_n
\le
\theta\sum_{i=1}^{n-1}\Delta_i.
\end{equation*}

\begin{lemma}[Presence of special indices]
\label{lem:special-index}
Every infinite $K$-special chain, with $K\ge L$, has infinitely many special indices.
\end{lemma}

\begin{proof}
Suppose otherwise. Then for all sufficiently large $n$,
$\Delta_n>\theta\sum_{i=1}^{n-1}\Delta_i$.
Set
\begin{equation*}
D_n:=\sum_{i=1}^n \Delta_i.
\end{equation*}
For all sufficiently large $n$, $D_n>(1+\theta)D_{n-1}$.
Thus there exists $c>0$ such that
$D_n\ge c(1+\theta)^n$
for all large $n$. On the other hand, \eqref{eq:special-deficit-upper} yields $D_n\le C(e_0)q^n.$
This is a contradiction since $1+\theta>q$. Hence we have proved the claim.
\end{proof}

We next exclude infinite special chains.

\begin{proposition}[No infinite special chains]
\label{prop:no-special-chains}
There exists $K_0<\infty$ such that, for every $K\ge K_0$, almost
surely there is no infinite $K$-special chain.
\end{proposition}

\begin{proof}
Fix a deterministic label $v_0$ and put $e_0=e_{v_0}$. For
$n\ge2$, let $\mathcal B_n(v_0)$ be the event that $e_0$ is the
initial edge of a $K$-special chain
$e_0\prec e_1\prec\cdots\prec e_n$ whose terminal index $n$ is a special index.

Suppose that such a chain occurs. Choose integers
$j_1,\ldots,j_n\ge K$ satisfying
\begin{equation*}
j_i\le \Delta(e_i)<j_i+1.
\end{equation*}
Since $n$ is a special index,
\begin{align}
j_n
\le \Delta(e_n)
\le
\theta\sum_{i=1}^{n-1}\Delta(e_i)<
\theta\sum_{i=1}^{n-1}j_i+\theta(n-1).
\end{align}
On the event that $e_0$ is $K$-special,
Lemma \ref{lem:special-chain-estimate} and Markov's inequality give
\begin{align}
&\P\!\left(
\mathcal B_n(v_0)\mid\mathscr P_{v_0}
\right)
\nonumber\\
&\le
C_L^n\mu^{-\Delta(e_0)/q}
\sum_{j_1,\ldots,j_{n-1}\ge K}
\mu^{-\frac{Q-1}{m}\sum_{i=1}^{n-1}j_i}
\sum_{\substack{j_n\ge K:\\
j_n<
\theta\sum_{i=1}^{n-1}j_i+\theta(n-1)}}
\mu^{j_n/m}.
\label{eq:special-prob-first}
\end{align}
There exists $C_m<\infty$ such that
$\sum_{j=K}^{J}\mu^{j/m}
\le
C_m\mu^{J/m},
\; J\ge K.$ Hence, using \eqref{eq:special-prob-first}, we obtain
\begin{align}
&\P\!\left(
\mathcal B_n(v_0)\mid\mathscr P_{v_0}
\right)
\le
C_mC_L^n
\mu^{-\Delta(e_0)/q}
\mu^{\theta(n-1)/m}
\sum_{j_1,\ldots,j_{n-1}\ge K}
\mu^{-\frac{Q-1-\theta}{m}
\sum_{i=1}^{n-1}j_i}.
\label{eq:special-prob-second}
\end{align}

Set
\begin{equation*}
\delta:=
\frac{Q-1-\theta}{m}>0\hspace{.3in}
\text{ and }\hspace{.3in}
B_K:=
\sum_{j=K}^{\infty}\mu^{-\delta j}.
\end{equation*}
The sum in \eqref{eq:special-prob-second} equals
$B_K^{n-1}$ and therefore, for a constant $C_1<\infty$,
\begin{equation}
\P\!\left(
\mathcal B_n(v_0)\mid\mathscr P_{v_0}
\right)
\le
C\mu^{-\Delta(e_0)/q}C_1^nB_K^{n-1}.
\label{eq:special-prob-final}
\end{equation}
If $e_0$ is $K$-special, then $\Delta(e_0)\ge K$, and hence $\mu^{-\Delta(e_0)/q}\le\mu^{-K/q}.$ Taking expectations in \eqref{eq:special-prob-final}
therefore yields
\begin{equation*}
\P(\mathcal B_n(v_0))
\le
C\mu^{-K/q}C_1^nB_K^{n-1}.
\end{equation*}
Since $B_K\downarrow0$ as $K\to\infty$, choose $K_0\ge L$ so large
that $C_1B_{K_0}<1.$ For every $K\ge K_0$,
\begin{equation*}
\sum_{n=2}^{\infty}\P(\mathcal B_n(v_0))<\infty.
\end{equation*}
By Borel--Cantelli, almost surely only finitely many such events occur for chains starting from $e_0$.

By Lemma \ref{lem:special-index}, an infinite $K$-special chain
would have infinitely many special indices and therefore, almost surely,
no such infinite $K$-special chain starts from $e_0$. Taking the union over all deterministic labels $v_0$ proves that almost surely there is no infinite $K$-special chain anywhere in the Galton--Watson tree.
\end{proof}
\subsection{Bounded deficits at new maxima times}
We now return to surviving rays. The special-chain estimates did not require edges to be new maxima; the connection with sequences of new maxima times is supplied by Corollary~\ref{cor:new-maxima-horizons}.

\begin{proof}[Proof of Theorem~\ref{thm:bounded-deficit}]
Choose $K\ge K_0$, where $K_0$ is given by
Proposition \ref{prop:no-special-chains}. Let $\zeta$ be a surviving ray. By
Lemma \ref{lem:new-maxima-inequality}, $\zeta$ has infinitely many edges at successive new maxima times:
\begin{equation*}
e_1\prec e_2\prec\cdots,
\qquad
t_i:=|e_i|.
\end{equation*}
By Corollary \ref{cor:new-maxima-horizons}, almost surely, simultaneously
for every surviving ray, all sufficiently late consecutive new maximum
edges satisfy
\begin{equation}
t_{i+1}\le H(e_i)+L.
\label{eq:special-new-maxima-window}
\end{equation}
Suppose that $\zeta$ contains only finitely many edges at new maxima times
satisfying $\Delta(e)\le K.$ Then there exists $i_0$ such that, for every $i\ge i_0$, $\Delta(e_i)>K$ and \eqref{eq:special-new-maxima-window} holds. Thus, for every $i\ge i_0$, the edge $e_{i+1}$ is a descendant of
$e_i$ satisfying
\begin{equation*}
t_i<t_{i+1}\le H(e_i)+L
\quad\text{ and }\quad
\Delta(e_{i+1})>K.
\end{equation*}
Hence $e_{i+1}$ is a $K$-special child of $e_i$. Therefore
\begin{equation*}
e_{i_0}\prec e_{i_0+1}\prec e_{i_0+2}\prec\cdots
\end{equation*}
is an infinite $K$-special chain. This contradicts Proposition \ref{prop:no-special-chains}. Consequently, every surviving ray contains infinitely many
edges at successive new maxima times $e$ such that
$\Delta(e)\le K.$ By \eqref{eq:deficit}, this is equivalent to $H(e)\ge q|e|-K.$
\end{proof}

\subsection{From bounded deficits to the burst speed}
\label{sec:accelerator}

We now show how the bounded-deficit times of Theorem~\ref{thm:bounded-deficit} can be used in an embedded branching construction. We first define the initial condition used in this construction.

\begin{definition}[Proper start]
\label{def:proper-start}
Fix $D>0$. Given a generation $s_0$, put
$s_1:=\lceil q s_0+D\rceil.$ An alive vertex $v_0$ at generation $s_0$ is called a \emph{proper start} for the parameters $D$ if
\begin{equation*}
S(v_0)\ge b_{s_1}+b_{s_0-1}.
\end{equation*}
\end{definition}

\subsubsection{Obtaining a proper start from a bounded deficit}

We show that a sufficiently late edge at a new maximum time with bounded deficit gives a proper start with probability bounded below independently of its generation.

\begin{lemma}[Obtaining a proper start]
\label{lem:proper-start}
Fix $D>0$ and let $K$ be the constant in Theorem \ref{thm:bounded-deficit}.  There exist $t_D<\infty$ and $p_D>0$ such that the
following holds.  Suppose that $v$ is an alive vertex at generation
$t\ge t_D$ and that its incoming edge $e$ satisfies
$H(e)\ge qt-K.$ Write
\begin{equation*}
s:=\lfloor qt-K\rfloor,
\qquad
s_1:=\lceil qs+D\rceil.
\end{equation*}
Then, conditionally on $\mathcal F_t$, with probability at least
$p_D$ there exists an alive descendant $w$ of $v$ at generation $s$
satisfying
\begin{equation*}
S(w)\ge b_{s_1}+b_{s-1}.
\end{equation*}
\end{lemma}

\begin{proof}
Since $s\le H(e)$ implies $X_e\ge b_s$ and the parent of $v$ is alive, we have
$S(v)\ge b_s+b_{t-1}.$ Set $d:=s-t$. Starting from $v$, we now retain a descendant path to generation $s-1$ if
every displacement along that path is at least $-b_{t-1}/d,$ which makes sure that the paths remain alive.
Write
\begin{equation*}
\varepsilon_t:=
\P\left(
X<-\frac{b_{t-1}}{d}
\right).
\end{equation*}
Since $d\le Ct$ for a deterministic constant $C$, Assumption
\ref{ass:Lstar} yields
$d\varepsilon_t=o(1).$

Since the retained 
genealogy is a Galton--Watson process with offspring mean
\begin{equation*}
\lambda_t=\mu(1-\varepsilon_t),
\qquad
\nu_t
:=
\E[\xi(\xi-1)](1-\varepsilon_t)^2
\le \nu,
\qquad
\nu:=\E[\xi(\xi-1)]<\infty,
\end{equation*}
where finiteness follows from \eqref{eq:GW2}.  Moreover, for all sufficiently large $t$,
$\lambda_t\ge \lambda_*>1.$

Let $W$ be the number of retained descendants at generation $s-1$ and put
$\ell:=s-t-1.$ Then $\E W=\lambda_t^\ell$ and 
by a standard computation using most recent common ancestor, we obtain
\begin{equation*}
\E W^2
=
\lambda_t^\ell
+
\nu_t
\sum_{r=0}^{\ell-1}
\lambda_t^{2\ell-r-2}.
\end{equation*}
Consequently,
\begin{align*}
\frac{\E W^2}{(\E W)^2}
&=
\lambda_t^{-\ell}
+
\frac{\nu_t}{\lambda_t^2}
\sum_{r=0}^{\ell-1}\lambda_t^{-r}
\le
1+
\frac{\nu}{\lambda_*^2}
\sum_{r=0}^{\infty}\lambda_*^{-r}
=:C<\infty
\end{align*}
uniformly for all sufficiently large $t$.

From each such retained particle at generation $s-1$, consider a child at
generation $s$ when its incoming displacement is at least $b_{s_1}$.
Let $Y$ be the number of such children.  By Assumption
\ref{ass:u-l-bound},
\begin{equation*}
\pi_t:=\P(X\ge b_{s_1})
\ge c_*\mu^{-s_1/m}
\end{equation*}
for all large $t$.  Consequently, using $1-q/m=1/q$, we have
\begin{align*}
\log_\mu \E Y
=
d-\frac{s_1}{m}+O(1)  =
\frac{s}{q}-t-\frac{D}{m}+O(1).
\end{align*}
Since
$s=qt-K+O(1)$, the right-hand side is bounded below uniformly in
$t$.  Hence $\E Y\ge c_{D,K}>0.$

Conditionally on $W$, the variable $Y$ is the sum of $W$ iid marked
offspring variables $\eta$ satisfying
$\E\eta=\mu\pi_t,
\;\;
\E[\eta(\eta-1)]=\nu\pi_t^2.$ Consequently,
\begin{align*}
(\E Y)^2
=
\mu^2\pi_t^2(\E W)^2 \quad
\text{ and }\quad\E Y^2
=
(\mu\pi_t)^2\E[W(W-1)]
+
(\nu\pi_t^2+\mu\pi_t)\E W.
\end{align*}
Hence, we obtain $\E Y^2/(\E Y)^2\le C_{D,K}$ and the
Paley--Zygmund inequality implies
$\P(Y>0\mid\mathcal F_t)\ge p_D>0.$ 
\end{proof}

\subsubsection{From a proper start to the strict regime $(q, Q)$}

We now use the proper start to accelerate from the critical value $q$
to the strict regime $(q, Q)$.  

\begin{proposition}[Accelerating trees]
\label{prop:accelerator}
Assume \eqref{eq:GW1}, \eqref{eq:GW2}, Assumptions
\ref{ass:Lstar} and \ref{ass:u-l-bound}, and let $m>4$. Fix
\begin{equation*}
q<\beta<R<Q.
\end{equation*}
Then there exist $D>0$, $s_*<\infty$ and $\rho>0$ such that the
following holds.

For every $s_0\ge s_*$, define recursively
\begin{equation}
s_{j+1}:=\left\lceil q s_j+DR^j\right\rceil,
\qquad j\ge0.
\label{eq:accelerator-times}
\end{equation}
Suppose that $v_0$ is a proper start for $D$ at generation $s_0$, that is,
\begin{equation}
S(v_0)\ge b_{s_1}+b_{s_0-1}.
\label{eq:proper-start-accelerator}
\end{equation}
Then, conditionally on $\mathcal F_{s_0}$, with probability at least $\rho$,
there exists a surviving ray having a sequence
$v_0\prec v_1\prec v_2\prec\cdots,$ with $|v_j|=s_j,$ such that
\begin{equation*}
S(v_j)\ge b_{s_{j+1}}+b_{s_j-1},
\qquad j\ge0.
\end{equation*}
Moreover, $s_{j+1}\ge\beta s_j$ for all sufficiently large $j$.
\end{proposition}
\begin{proof}
We construct a tree $\T_{\mathrm{acc}}$ embedded in the original
Galton--Watson tree. The root of $\T_{\mathrm{acc}}$ is $v_0$, and
the vertices at level $j$ of $\T_{\mathrm{acc}}$ are vertices of the
original tree at generation $s_j$. We call such a vertex $v$ active if
$S(v)\ge b_{s_{j+1}}+b_{s_j-1}.$
Thus, $v_0$ is active by \eqref{eq:proper-start-accelerator}.

Set $d_j:=s_{j+1}-s_j.$ Starting from an active vertex $v$ at generation $s_j$, we retain a
descendant path
\begin{equation*}
v=u_0\prec u_1\prec\cdots\prec u_{d_j-1},
\qquad |u_r|=s_j+r,
\end{equation*}
whenever
\begin{equation}
X_{u_r}\ge
-\frac{b_{s_j-1}}{d_j},
\qquad 1\le r<d_j.
\label{eq:stable-edge-accelerator}
\end{equation}
From every such retained vertex $u_{d_j-1}$ at generation $s_{j+1}-1$,
we retain a child $w$ at generation $s_{j+1}$ if
\begin{equation}
X_w\ge b_{s_{j+2}}.
\label{eq:terminal-jump-accelerator}
\end{equation}
The vertices $w$ obtained are the children of $v$ in
$\T_{\mathrm{acc}}$.

We first observe that every vertex retained in this construction is
alive. Indeed, suppose $v$ is active and let
$v=u_0\prec u_1\prec\cdots\prec u_{d_j-1}$ be a retained path. For every $1\le r<d_j$,
\begin{align*}
S(u_r)
&=
S(v)+\sum_{\ell=1}^r X_{u_\ell}\ge
b_{s_{j+1}}+b_{s_j-1}
-r\frac{b_{s_j-1}}{d_j}\ge b_{s_{j+1}}.
\end{align*}
Since $|u_r|<s_{j+1}$ and the barrier is increasing, $S(u_r)\ge b_{s_{j+1}}\ge b_{|u_r|}.$
If $w$ further satisfies \eqref{eq:terminal-jump-accelerator}, then
\begin{equation*}
S(w)
\ge b_{s_{j+1}}+b_{s_{j+2}}
\ge b_{s_{j+1}-1}+b_{s_{j+2}},
\end{equation*}
and hence $w$ is active at level $j+1$. Thus, survival of
$\T_{\mathrm{acc}}$ implies survival of the killed branching random
walk.

We now study the survival of $\T_{\mathrm{acc}}$. Let $Y_j$ denote the number of children in $\T_{\mathrm{acc}}$
produced by a single active vertex at level $j$.
We state the required one-block estimate and prove it after completing the argument.

\begin{lemma}
\label{lem:accelerator-one-block}
There exist constants $D_0<\infty$, $s_*<\infty$ and $C<\infty$
such that, for every fixed $D\ge D_0$, taking $s_*$ large enough, the following holds uniformly for $s_0\ge s_*$ and
$j\ge0$. If $m_j:=\E Y_j,$ then
\begin{equation*}
\log_\mu m_j
\ge
\frac{Q-R}{m}DR^j-C,
\end{equation*}
and
\begin{equation}
\E Y_j^2\le C m_j^2.
\label{eq:good-offspring-second}
\end{equation}
In particular, $D$ may be chosen sufficiently large so that
$m_j\ge2$ for $j\ge0.$
\end{lemma}

Returning to the proof, fix such a $D$ and choose $s_*$ accordingly. Let $Z_j$ denote the number of active vertices
at level $j$ of $\T_{\mathrm{acc}}$, so that $Z_0=1$. By
branching property, conditionally on
$\mathcal F_{s_0}$, we have
\begin{align}\label{eq:good-mean-lower-j}
M_j:=\E[Z_j\mid\mathcal F_{s_0}]\geq 2^{j}.
\end{align}
Furthermore, conditionally on $Z_j$,
$\E[Z_{j+1}^2\mid Z_j,\mathcal F_{s_0}]
=
m_j^2Z_j^2
+
\operatorname{Var}(Y_j)Z_j.$ Since $M_{j+1}^2=m_j^2M_j^2$, we have the following recursion by \eqref{eq:good-offspring-second}:
\begin{align*}
\frac{\E[Z_{j+1}^2\mid\mathcal F_{s_0}]}{M_{j+1}^2}
&=
\frac{\E[Z_j^2\mid\mathcal F_{s_0}]}{M_j^2}
+
\frac{\operatorname{Var}(Y_j)}{m_j^2M_j}\le
\frac{\E[Z_j^2\mid\mathcal F_{s_0}]}{M_j^2}
+
\frac{C}{M_j}.
\end{align*}
Iterating this
inequality and using \eqref{eq:good-mean-lower-j}, we get
\begin{align*}
\frac{\E[Z_j^2\mid\mathcal F_{s_0}]}{M_j^2}
&\le
1+C\sum_{i=0}^{j-1}\frac1{M_i}
\le
1+2C.
\end{align*}
Therefore the Paley--Zygmund inequality along with the fact that events $\{Z_j>0\}$ are decreasing in $j$, yields
\begin{align*}
\P(Z_j>0\text{ for every }j\mid\mathcal F_{s_0})
&=\lim_{j\to\infty}\P(Z_j>0\mid\mathcal F_{s_0})\ge\frac{1}{1+2C}.
\end{align*}
On this event the embedded tree has vertices at every level and is locally finite, so we get an active ray. Finally, set $u_j:=s_j/R^j$.  From \eqref{eq:accelerator-times}, 
\begin{align*}
u_{j+1}=\frac qR u_j+\frac DR+O(R^{-j-1}).
\end{align*}
Since $q/R<1$, it follows that $u_j\to D/(R-q)$, and therefore $s_{j+1}/s_j\to R$.  Since $\beta<R$, we also have $s_{j+1}\ge\beta s_j$ for all sufficiently large $j$.  This completes the proof of the proposition.
\end{proof}
\begin{proof}[Proof of Lemma~\ref{lem:accelerator-one-block}]
Fix an active vertex at generation $s_j$ and write
\begin{align*}
d_j:=s_{j+1}-s_j,
\qquad
\varepsilon_j
:=
\P\left(
X<-\frac{b_{s_j-1}}{d_j}
\right).
\end{align*}
Since $R>q$ and $s_0$ is sufficiently large, there is a constant
$C_R<\infty$ such that $d_j\le C_Rs_j,$
 $j\ge0$.

Assumption \ref{ass:Lstar} therefore yields
$d_j\varepsilon_j=o(1)$
uniformly in $j$ as $s_0\to\infty$.

Let $W_j$ be the number of descendants at generation $s_{j+1}-1$
whose intervening displacements satisfy
\eqref{eq:stable-edge-accelerator}. As in the proof of
Lemma \ref{lem:proper-start}, these descendants form a
Galton--Watson process with offspring mean
$\lambda_j:=\mu(1-\varepsilon_j).$
Hence
\begin{align*}
\E W_j=\lambda_j^{d_j-1}.
\end{align*}
Moreover, since $\E\xi^2<\infty$ and $\lambda_j$ is uniformly bounded
away from $1$ for all sufficiently large $s_0$, the same
second-moment calculation as in Lemma \ref{lem:proper-start} yields
\begin{equation}
\frac{\E W_j^2}{(\E W_j)^2}\le C_0
\label{eq:accelerator-W-second}
\end{equation}
uniformly in $j$ and in all sufficiently large $s_0$.

Now write
$\pi_j
:=
\P(X\ge b_{s_{j+2}}).$
Each child of a retained particle satisfies the terminal-jump
condition independently with probability $\pi_j$. Consequently,
\begin{equation*}
m_j:=\E Y_j
=
\mu\pi_j\E W_j
=
\mu\pi_j\lambda_j^{d_j-1}.
\end{equation*}
By Assumption \ref{ass:u-l-bound},
$\pi_j\asymp\mu^{-s_{j+2}/m}.$
Together with $d_j\varepsilon_j=o(1)$, this yields
\begin{equation}
\log_\mu m_j
\ge
d_j-\frac{s_{j+2}}m-C
\label{eq:mj-bound}
\end{equation}
for some constant $C<\infty$. Using
$1-q/m=1/q=Q/m,$ we obtain
\begin{align*}
s_{j+1}-s_j-\frac{s_{j+2}}m&=
\left(1-\frac qm\right)s_{j+1}
-s_j-\frac{DR^{j+1}}m+O(1)\\
&=
\frac{DR^j}{q}
-\frac{DR^{j+1}}m+O(1)=
\frac{Q-R}{m}DR^j+O(1).
\end{align*}
Substituting this into \eqref{eq:mj-bound} yields
\begin{equation}
\log_\mu m_j
\ge
\frac{Q-R}{m}DR^j-C.
\end{equation}
Since $R<Q$, choosing $D$ sufficiently large yields
$m_j\ge2,
\; j\ge0.$
It remains to establish the second-moment bound for $Y_j$. Put
$\nu:=\E[\xi(\xi-1)]<\infty,$ using branching property and terminal displacements yields,
conditionally on $W_j$,
\begin{align*}
\E[Y_j\mid W_j]
&=
\mu\pi_j W_j,\\
\E[Y_j(Y_j-1)\mid W_j]
&=
\nu\pi_j^2W_j
+
\mu^2\pi_j^2W_j(W_j-1).
\end{align*}
Since $m_j=\mu\pi_j\E W_j$, it follows that
\begin{align*}
\frac{\E Y_j^2}{m_j^2}
&\le
\frac{\E W_j^2}{(\E W_j)^2}
+
\frac{\nu}{\mu^2\E W_j}
+
\frac1{m_j}.
\end{align*}
By \eqref{eq:accelerator-W-second}, the first term is uniformly
bounded. Moreover, for all sufficiently large $s_0$,
$\E W_j\ge1$, while $m_j\ge2$. Hence $\E Y_j^2\le C m_j^2$
uniformly in $j$ and in all sufficiently large $s_0$. This proves
the lemma.
\end{proof}

\subsubsection{Proof of the burst speed}
The upper bound follows directly from the alive cloud. Indeed, for every
surviving ray $\zeta$ and every $n$,
$S_n(\zeta)\le M_n.$
By the cloud upper bound proved in Section~\ref{sec:cloud-upper},
\begin{equation*}
\limsup_{n\to\infty}
\frac{\ainv(S_n(\zeta))}{n}
\le
\gamma(m)
\qquad
\text{for every }\zeta\in\Rinf
\end{equation*}
almost surely on $\cS$.

\begin{proof}[Proof of the burst lower bound in Theorem~\ref{thm:main}]
Fix $q<\beta<R<Q.$ Let $D$, $s_*$ and $\rho$ be given by Proposition
\ref{prop:accelerator}, and apply Lemma \ref{lem:proper-start} with
this value of $D$.

There exists $p_\beta>0$ such that, for all sufficiently large $t$,
if $\mathcal F_t$ contains an alive vertex whose incoming edge $e$
satisfies
\begin{equation}
H(e)\ge qt-K,
\label{eq:ray-new-maxima-1}
\end{equation}
then, conditionally on $\mathcal F_t$, with probability at least
$p_\beta$ there is an infinite ray $\zeta$ for which
\begin{equation*}
S_{s_j}\ge b_{s_{j+1}}
=
a(s_{j+1}/m)
\end{equation*}
and $s_{j+1}\ge\beta s_j$ for all sufficiently large $j$. Hence
$\ainv(S_{s_j})/s_j
\ge
\beta/m$ for all sufficiently large $j$, and therefore
\begin{equation*}
\limsup_{n\to\infty}
\frac{\ainv(S_n(\zeta))}{n}
\ge
\frac{\beta}{m}.
\end{equation*}

Set
\begin{equation*}
A_\beta
:=
\left\{
\sup_{\zeta\in\Rinf}
\limsup_{n\to\infty}
\frac{\ainv(S_n(\zeta))}{n}
\ge
\frac{\beta}{m}
\right\}.
\end{equation*}
Then $Z_t^{(\beta)}
:=
\P(A_\beta\mid\mathcal F_t)$ is a bounded martingale and therefore
\begin{equation*}
\lim_{t\to\infty}Z_t^{(\beta)}=\1_{A_\beta}
\qquad\text{almost surely}.
\end{equation*}
On $\cS$, fix a surviving ray. By Theorem
\ref{thm:bounded-deficit}, this ray contains infinitely many new maximum
edges satisfying \eqref{eq:ray-new-maxima-1}. At each such
generation,
\begin{equation*}
Z_t^{(\beta)}\ge p_\beta.
\end{equation*}
Thus $Z_t^{(\beta)}$ is bounded away from zero along an unbounded
sequence of generations. Consequently $A_\beta^c$ cannot occur on
$\cS$, and therefore $\P(A_\beta\mid\cS)=1.$ Now choose $\beta_j\uparrow Q$ with $q<\beta_j<Q$. Almost surely on
$\cS$, all the events $A_{\beta_j}$ occur. Hence
\begin{align*}
\sup_{\zeta\in\Rinf}
\limsup_{n\to\infty}
\frac{\ainv(S_n(\zeta))}{n}
&\ge
\sup_{j\ge1}\frac{\beta_j}{m}=
\frac{Q}{m}
=
\gamma(m).
\end{align*}
The reverse inequality follows from the cloud upper bound and \eqref{eq:speed-order}.
\end{proof}

\subsection{The sustainable speed}
\label{sec:sustained}

Throughout this section, for $r>0$, write
$b_n^{(r)}:=a(n/r)$ and retain $b_n=b_n^{(m)}$. The following lemma gives a sufficient initial height for constructing a ray above the stronger $r$-barrier.

\begin{lemma}[Restart above an $r$-barrier]
\label{lem:r-restart}
Assume \eqref{eq:GW1}, \eqref{eq:GW2}, Assumptions~\ref{ass:Lstar} and \ref{ass:u-l-bound}, and fix $4<r<m$. There exist $s_r<\infty$ and $\rho_r>0$ such that, if $v$ is an
alive vertex at generation $s\ge s_r$ satisfying
\begin{equation*}
S(v)\ge b_{2s}^{(r)}+b_{s-1},
\end{equation*}
then, conditionally on $\mathcal F_s$, with probability at least
$\rho_r$ there is an infinite descendant ray $\zeta$ satisfying
\begin{equation*}
S_n(\zeta)\ge b_n^{(r)},
\qquad n\ge s.
\end{equation*}
\end{lemma}

\begin{proof}
Set $t_0=s$ and $t_{k+1}=2t_k$. At generation $t_k$, call a particle
active if
\begin{equation*}
S\ge b_{t_{k+1}}^{(r)}+b_{t_k-1}.
\end{equation*}
Starting from an active particle, keep descendant paths to generation
$t_{k+1}-1$ whose increments are all at least
$-b_{t_k-1}/t_k$,
and require the final displacement at $t_{k+1}$ to be at least
$b_{t_{k+2}}^{(r)}$.

Let $Y_k$ be the number of active children produced by one active
particle. Assumption \ref{ass:Lstar} and the same second-moment
argument as in Lemma \ref{lem:accelerator-one-block} give
\begin{align*}
\log_\mu\E Y_k
&\ge
t_k-\frac{t_{k+2}}r-C=
\left(1-\frac4r\right)t_k-C,
\end{align*}
and
$\E Y_k^2\le C_r(\E Y_k)^2$.

Since $r>4$, for $s$ sufficiently large we have $\E Y_k\ge2$ for
every $k$. Hence the embedded process survives with probability
bounded below by some $\rho_r>0$, by similar arguments as in Proposition \ref{prop:accelerator}. Moreover, during the block $[t_k,t_{k+1})$,
\begin{equation*}
S_n\ge b_{t_{k+1}}^{(r)}\ge b_n^{(r)},
\end{equation*}
and a successful final jump produces an active particle at the next
level. The resulting infinite ray therefore stays above
$b_n^{(r)}$ for all $n\ge s$.
\end{proof}

\begin{proof}[Proof of the sustainable-speed in Theorem~\ref{thm:main}]
The upper bound uses the extinction result from \citet{SamantC2026}. For every fixed $r<4$, almost surely there is no infinite ray whose position remains above $b_n^{(r)}$ for all sufficiently large $n$. Taking the intersection over rational $r\in(0,4)$ gives, simultaneously for all surviving rays $\zeta$,
\begin{align*}
\liminf_{n\to\infty}\frac{\ainv(S_n(\zeta))}{n}\le\frac14.
\end{align*}
For the reverse inequality, fix $4<r<m$. Since $r>4$, $2m/r<Q.$ Choose
\begin{equation*}
\max\left\{q,\frac{2m}{r}\right\}
<
\beta<R<Q.
\end{equation*}
By Lemma \ref{lem:proper-start} and Proposition
\ref{prop:accelerator}, every sufficiently late bounded-deficit edge
$H(e)\ge q|e|-K$ has a probability bounded away from zero of producing generation $s_j$ satisfying
\begin{equation*}
S_{s_j}\ge b_{s_{j+1}}+b_{s_j-1},
\qquad
s_{j+1}\ge\beta s_j
\end{equation*}
for all sufficiently large $j$. Since $\beta/m>2/r$,
\begin{equation*}
S_{s_j}
\ge
a(\beta s_j/m)+b_{s_j-1}
\ge
b_{2s_j}^{(r)}+b_{s_j-1}.
\end{equation*}
Lemma \ref{lem:r-restart} therefore implies that there exists $p_r>0$ such that every sufficiently late
bounded-deficit new maximum has conditional probability at
least $p_r$ of producing a ray satisfying
$S_n(\zeta)\ge b_n^{(r)}
$
for all sufficiently large $n$. Let
\begin{equation*}
A_r
:=
\left\{
\exists\,\zeta\in\Rinf:
S_n(\zeta)\ge b_n^{(r)}
\text{ for all sufficiently large }n
\right\}.
\end{equation*}
Then,
$Z_t^{(r)}:=\P(A_r\mid\mathcal F_t)
\longrightarrow
\mathbf 1_{A_r}$
almost surely. On $\cS$, Theorem~\ref{thm:bounded-deficit} gives infinitely
many bounded-deficit edges at new maxima times satisfying $H(e)\ge q|e|-K.$
Hence $Z_t^{(r)}\ge p_r$ along an unbounded sequence of generations. Therefore
$\P(A_r\mid\cS)=1.$
Consequently,
\begin{equation*}
\sup_{\zeta\in\Rinf}
\liminf_{n\to\infty}
\frac{\ainv(S_n(\zeta))}{n}
\ge
\frac14
\end{equation*}
almost surely on $\cS$, after taking a countable sequence $r_k\downarrow4$ with $4<r_k<m$ and using the inequality $s_{\rm sust}\ge1/r_k$ on $A_{r_k}$.  Together with the upper bound, this proves \eqref{eq:main-equality}.
\end{proof}

\section{Cloud Speed and Growth of the Alive Cloud}\label{sec:cloud-convergence}

We prove Theorem~\ref{thm:cloud-convergence} and determine the exponential
growth rate of the alive population.
The cloud upper bound is already contained in Theorem~\ref{thm:main}. For the lower bounds we use Lemma~\ref{lem:scale-separation}
and construct a subtree with controlled negative displacements.

For $K>0$, retain an edge $e$ when $X_e\ge-K$.
The retained genealogy is a Galton--Watson tree
with mean offspring
$$
\mu_K=\mu p_K,\qquad p_K:=\P(X\ge-K),
$$
and $\mu_K\uparrow\mu$ as $K\to\infty$.
Let $Z_h^{(K)}$ denote its generation-$h$ population started from one particle.

Suppose $\mu_K>1$. By \eqref{eq:GW2} and the
standard martingale convergence theorem for
Galton--Watson processes
\citep[Chapter~I]{AthreyaNey1972},
$
W_h:=Z_h^{(K)}/\mu_K^h
$
is a nonnegative martingale which converges almost surely and in
$L^2$ to a random variable $W_\infty$ with
$\E W_\infty=1$. In particular,
$\P(W_\infty>0)>0$.
Fix $\sigma<\log_\mu\mu_K$. On
$\{W_\infty>0\}$, almost surely, for all
sufficiently large $h$,
$$
Z_h^{(K)}
\ge \frac{W_\infty}{2}\mu_K^h
\ge \mu^{\sigma h}.$$
Consequently,
$\{W_\infty>0\}
\subseteq
\bigcup_{h_0\ge1}
\left\{
Z_h^{(K)}\ge\mu^{\sigma h}
\text{ for every }h\ge h_0
\right\}.$ 
Since $\P(W_\infty>0)>0$, there exist
deterministic constants $h_0<\infty$
and $p_*>0$ such that
\begin{equation}
\P\left(
Z_h^{(K)}\ge\mu^{\sigma h}
\text{ for every }h\ge h_0
\right)\ge p_*.
\label{eq:retained-population-growth}
\end{equation}
For $\beta>0$, put
$$
H_t(\beta):=\{\exists v\in\mathcal A_t:S(v)\ge a(\beta t)\},
\qquad
N_t(\beta):=\#\{v\in\mathcal A_t:S(v)\ge a(\beta t)\}.
$$
Fix
\begin{equation}
1-\gamma(m)<\beta<\gamma(m),\qquad
\beta<\vartheta<1,\qquad
\frac{\vartheta}{\vartheta-\beta}<A<B<m\beta.
\label{eq:propagation-params}
\end{equation}
Such choices exist since $m\beta(1-\beta)>1$. Let
\begin{equation*}
G_{\beta,\vartheta,B}:=
\left\{H_n(\beta)\text{ eventually},\quad
\liminf_{n\to\infty}\frac{\log_\mu(A_n\vee1)}{n}
\ge\vartheta\left(1-\frac1B\right)\right\}.
\end{equation*}

\begin{lemma}
\label{lem:high-particle-propagation}
Assume \eqref{eq:GW1}, \eqref{eq:GW2}, constraints \eqref{C1}--\eqref{C3}, and $m>4$. For the parameters in \eqref{eq:propagation-params}, there exist $t_0<\infty$ and $p>0$ such that, for all $t\ge t_0$,
\begin{equation*}
\P(G_{\beta,\vartheta,B}\mid\cF_t)
\ge p\,\mathbf 1_{H_t(\beta)}.
\end{equation*}
\end{lemma}

\begin{proof}
Choose $K$ such that
$
\vartheta<\log_\mu\mu_K<1.
$
We first show that a $\beta$-high particle at sufficiently late generation $u$, produces many $\beta$-high descendants at every generation in $[Au,Bu]$.

Let $v\in\mathcal A_u$ satisfy $S(v)\ge a(\beta u)$. By Lemma~A.1, for all sufficiently large $u$, every descendant $w$ of $v$ in the retained genealogy, with $|w|\le\lfloor Bu\rfloor$, is alive. Indeed, along the retained paths,
$$
S(w)\ge a(\beta u)-K(|w|-u)
\ge b_{\lceil Bu\rceil}\ge b_{|w|},
$$
and the same inequality holds at the intermediate generations.
\medskip

For an integer $j\in[Au,Bu]$, let $Y_{v,j}$ denotes the number of retained descendants $w$ of $v$ at generation $j$ satisfying $X_w\ge a(\beta j)$. Conditionally on the retained genealogy below $v$, the displacements on its edges are independent with law $\mathcal L(X\mid X\ge-K)$. Consequently, for all sufficiently large $j$,
$$
Y_{v,j}\ \sim\ \operatorname{Bin}\left(
Z_{j-u}^{(K)},\frac{\Fbar(a(\beta j))}{p_K}\right)
$$
under this conditional law. By \eqref{C2}, for every $\varepsilon>0$ and all sufficiently large $j$,
$
\Fbar(a(\beta j))\ge\mu^{-(\beta+\varepsilon)j}.
$
Since $\vartheta(1-1/A)>\beta$, we may choose $\varepsilon,\delta>0$ so small that
$\nu\left(1-1/A\right)-\beta-\varepsilon>2\delta$.

On the event in \eqref{eq:retained-population-growth}, for all sufficiently large $u$ and every integer $j\in[Au,Bu]$,
\begin{align*}
\E[Y_{v,j}\mid\text{retained genealogy}]
&=Z_{j-u}^{(K)}\frac{F(a(\beta j))}{p_K}\ge\mu^{\nu(j-u)-(\beta+\varepsilon)j}
\ge\mu^{2\delta j}.
\end{align*}
A Chernoff bound for a binomial random variable then gives
$$
\P\left(Y_{v,j}<\mu^{\delta j}\mid
\text{retained genealogy}\right)
\le\exp\{-c\mu^{2\delta j}\}
$$
for some $c>0$, uniformly over retained genealogies satisfying \eqref{eq:retained-population-growth}. Summing over $j\in[Au,Bu]$ and using \eqref{eq:retained-population-growth}, we find $u_0<\infty$ and $p_0>0$ such that, for every $u\ge u_0$ and every $\beta$-high alive vertex $v$ at generation $u$,
\begin{equation}
\P\left(
\begin{array}{c}
\text{$v$ has at least $\mu^{\delta j}$ alive, $\beta$-high}\\
\text{descendants at every integer $j\in[Au,Bu]$}
\end{array}
\,\middle|\,\cF_u\right)\ge p_0.
\label{eq:cloud-one-seed}
\end{equation}
We next use the number of high particles to continue this construction. Suppose that $N_u(\beta)\ge\mu^{\delta u}$ for some $u\ge u_0$. Conditional on $\cF_u$, by branching property at generation $u$ and \eqref{eq:cloud-one-seed}, the probability that none of these particles produces the required high descendants throughout $[Au,Bu]$ is at most
$
(1-p_0)^{N_u(\beta)}\le\exp\{-p_0\mu^{\delta u}\}.
$
Thus, if
$$
E_u:=\left\{
N_u(\beta)\ge\mu^{\delta u},\quad
N_j(\beta)<\mu^{\delta j}
\text{ for some }j\in[Au,Bu]\cap\mathbb N
\right\},
$$
then
$\P(E_u\mid\cF_u)\le\exp\{-p_0\mu^{\delta u}\}$.
Hence, for every $t\ge u_0$,
\begin{equation}
\P\left(\bigcup_{u\ge t}E_u\,\middle|\,\cF_t\right)
\le\sum_{u\ge t}\exp\{-p_0\mu^{\delta u}\}
\rightarrow 0.
\label{eq:cloud-summable-failures}
\end{equation}

On $H_t(\beta)$, choose a $\beta$-high alive vertex
$v$ at generation $t$. By \eqref{eq:cloud-one-seed},
conditionally on $\cF_t$, with probability at least
$p_0$,
\begin{equation}
N_j(\beta)\ge\mu^{\delta j}
\qquad
\text{for every integer }j\in[At,Bt].
\label{eq:cloud-first-interval}
\end{equation}

We show that, if this event occurs and none of the
events $E_u$, $u\ge t$, occur, then the same lower bound holds at all sufficiently large generations.
Indeed, by the definition of $E_u$, on $E_u^c$,
whenever $N_u(\beta)\ge\mu^{\delta u},$ we necessarily have
\begin{equation}
N_j(\beta)\ge\mu^{\delta j}
\qquad
\text{for every integer }j\in[Au,Bu].
\label{eq:restart}
\end{equation}
Put
$
L_1=\lceil At\rceil$, and 
$U_1=\lfloor Bt\rfloor.
$
By \eqref{eq:cloud-first-interval}, the required
lower bound holds throughout $[L_1,U_1]$.
We may therefore apply
\eqref{eq:restart} at each integer
generation $u\in[L_1,U_1]$. Since $B>A$, for all sufficiently large $u$,
$$
\lceil A(u+1)\rceil\le\lfloor Bu\rfloor+1.
$$
Their union contains every integer
in
$
[L_2,U_2]$, where
$L_2=\lceil AL_1\rceil$ and
$U_2=\lfloor BU_1\rfloor.
$
Consequently,
$$
N_j(\beta)\ge\mu^{\delta j}
\qquad
\text{for all integers }j\in[L_2,U_2].
$$
Repeating this argument, define
$L_{r+1}=\lceil AL_r\rceil,\;
U_{r+1}=\lfloor BU_r\rfloor,
\; r\ge1.$
By induction,
$
N_j(\beta)\ge\mu^{\delta j},
\;
\text{for every integer }j\in[L_r,U_r].
$
Moreover,
$$
L_r\le A^r\left(t+\frac1{A-1}\right),
\qquad
U_r\ge B^r\left(t-\frac1{B-1}\right).
$$
Since $B>A$, we have $U_r\ge L_{r+1}-1$
for all sufficiently large $r$.
The intervals therefore eventually overlap
and cover all sufficiently large integers and hence it follows that
\begin{equation}
N_n(\beta)\ge\mu^{\delta n}
\qquad\text{for all sufficiently large }n.
\label{eq:cloud-high-population}
\end{equation}

\medskip

It remains to obtain the population lower bound. Put $u_n=\lceil n/B\rceil$ and $h_n=n-u_n$. Since $\mu_K>\mu^\vartheta$, \eqref{eq:retained-population-growth} gives a constant $p_1>0$ such that, for all sufficiently large $n$,
$$
\P\left(Z_{h_n}^{(K)}\ge\mu^{\vartheta h_n}\right)\ge p_1.
$$
By Lemma~A.1, every retained descendant up to generation $n$ of a $\beta$-high vertex at generation $u_n$ is alive, since $n\le Bu_n$. Therefore, conditionally on $\cF_{u_n}$, using branching property, we have
\begin{align*}
&\P\left(A_n<\mu^{\vartheta h_n},\,
N_{u_n}(\beta)\ge\mu^{\delta u_n}\right)\le \exp\{-p_1\mu^{\delta u_n}\}.
\end{align*}
By Borel--Cantelli, almost surely on the event \eqref{eq:cloud-high-population},
$
A_n\ge\mu^{\vartheta(n-u_n)}$
for all sufficiently large $n$.
Since $u_n/n\to1/B$, this yields
$$
\liminf_{n\to\infty}\frac{\log_\mu(A_n\vee1)}{n}
\ge\vartheta\left(1-\frac1B\right).
$$

We have shown that, outside a null event, $G_{\beta,\vartheta,B}$ occurs on $H_t(\beta)$ whenever \eqref{eq:cloud-first-interval} holds and none of the events $E_u$, $u\ge t$, occurs. By \eqref{eq:cloud-one-seed} and \eqref{eq:cloud-summable-failures},
$$
\P(G_{\beta,\vartheta,B}\mid\cF_t)
\ge(p_0-\varepsilon_t)\mathbf1_{H_t(\beta)}.
$$
Here $\varepsilon_t$ denotes the right-hand side of \eqref{eq:cloud-summable-failures}. Taking $t_0$ sufficiently large that $\varepsilon_t\le p_0/2$ for $t\ge t_0$ proves the lemma with $p=p_0/2$. 
\end{proof}

\begin{proof}[Proof of Theorem~\ref{thm:cloud-convergence}]
Write $\gamma=\gamma(m)$ and fix
$1-\gamma<\beta<\gamma$. By Theorem~\ref{thm:main},
$H_t(\beta)$ occurs infinitely often on $\cS$.
For admissible $\vartheta,A,B$, put $Y_t:=\P(G_{\beta,\vartheta,B}\mid\cF_t).$ Since $G_{\beta,\vartheta,B}\in\cF_\infty$,
martingale convergence gives
$$
\lim_{t \rightarrow \infty}Y_t=\mathbf1_{G_{\beta,\vartheta,B}}
\qquad\text{almost surely}.
$$
By Lemma~\ref{lem:high-particle-propagation},
$Y_t\ge p$ on $H_t(\beta)$ for all sufficiently
large $t$. Consequently,
$$
\limsup_{t\to\infty}Y_t\ge p
\qquad\text{almost surely on }\cS.
$$
This implies that $G_{\beta,\vartheta,B}$ holds
almost surely on $\cS$. In particular,
$$
\liminf_{n\to\infty}
\frac{\ainv(M_n)}{n}\ge\beta
\qquad\text{almost surely on }\cS.
$$
Letting $\beta\uparrow\gamma$ through a
countable sequence and using the upper bound
in Theorem~\ref{thm:main} completes the proof.
\end{proof}

\begin{proof}[Proof of Theorem~\ref{thm:alive-population-growth}]
Write $\gamma=\gamma(m)$. We first prove the upper bound.
Fix $\varepsilon>0$ and choose
$\gamma<\beta<\gamma+\varepsilon$. Put
$
x_n=a(\beta(n+1)),\text{ and } p_n=\Fbar(x_n).
$
By \eqref{C2},
$
p_n=\mu^{-\beta(n+1)+o(n)}.
$
Let $D_n$ be the event that
$A_n\ge\mu^{(\gamma+\varepsilon)n}$ and no edge
from an alive particle at generation $n$ has
displacement at least $x_n$.

Put $p_+=\P(\xi\ge1)>0$. For each alive particle at generation $n$, the probability that it has at least one child whose displacement is at least
$x_n$ is at least $p_+p_n$. By the branching
property, these events are independent
conditionally on $\cF_n$. Hence
\begin{align*}
\P(D_n\mid\cF_n)
&\le
\mathbf1_{\{A_n\ge\mu^{(\gamma+\varepsilon)n}\}}
(1-p_+p_n)^{A_n}\le
\exp\left\{
-p_+\mu^{(\gamma+\varepsilon)n}p_n
\right\}.
\end{align*}
Since $\beta<\gamma+\varepsilon$, there exists
$c>0$ such that, for all sufficiently large $n$,
$
\P(D_n)\le\exp\{-\mu^{cn}\}.
$
In particular, $\sum_n\P(D_n)<\infty$.
Since $\beta>\gamma>1/m$, we have
$x_n>b_{n+1}$ for all sufficiently large $n$.
An edge with displacement at least $x_n$ from
an alive particle at generation $n$ therefore
produces an alive child with position at least
$x_n$. By the cloud upper bound in
Theorem~\ref{thm:main}, almost surely on $\cS$,
$
M_{n+1}<x_n~
\text{eventually}.
$
Thus, on $\cS$, whenever
$A_n\ge\mu^{(\gamma+\varepsilon)n}$ at a
sufficiently large generation, $D_n$ occurs.
By Borel--Cantelli, we have
$$
\limsup_{n\to\infty}
\frac{\log_\mu A_n}{n}
\le\gamma+\varepsilon
\qquad\text{on }\cS.
$$
Letting $\varepsilon\downarrow0$ proves the
upper bound. For the lower bound, for each admissible choice of
$\beta, \vartheta, A, B$, Lemma \ref{lem:high-particle-propagation}
gives
$$
\liminf_{n\to\infty}
\frac{\log_\mu A_n}{n}
\ge
\vartheta\left(1-\frac1B\right)
\qquad\text{almost surely on }\cS.
$$
Let $\beta\uparrow\gamma$, $\vartheta\uparrow1$
and $B\uparrow m\gamma$ through admissible
choices. Since $m\gamma(1-\gamma)=1$, we obtain
$$
\liminf_{n\to\infty}
\frac{\log_\mu A_n}{n}
\ge 1-\frac1{m\gamma}
=\gamma
\qquad\text{on }\cS.
$$
This proves the claim.
\end{proof}

\begin{proof}[Proof of Corollary~\ref{cor:unkilled-speeds}]
Fix $m>4$ for which the assumptions hold and write
$p_m:=\P(\cS)>0$. By the branching property, conditionally
on $\cF_n$, the probability that none of the subtrees
rooted at generation $n$ contains a surviving copy of the
killed branching random walk is $(1-p_m)^{Z_n}$.
Since $Z_n\to\infty$ almost surely on survival of the
Galton--Watson tree, almost surely on survival such a
subtree exists.

By Lemma~\ref{lem:scale-separation}, finite shifts in
position and generation do not affect the asymptotic
speeds. Hence, by Theorems~\ref{thm:main}
and~\ref{thm:cloud-convergence},
\begin{equation*}
    s_{\rm burst}\ge\gamma(m),\qquad
    s_{\rm sust}\ge\frac14,\qquad
    \liminf_{n\to\infty}\frac{\ainv(M_n^{(0)})}{n}
    \ge\gamma(m).
\end{equation*}
Letting $m\to\infty$ through a countable sequence gives
the required lower bounds.

For the upper bounds, by \eqref{C3}, for every $\beta>1$,
\begin{equation*}
    \P(M_n^{(0)}\ge a(\beta n))
    \le n\mu^n\Fbar\left(\frac{a(\beta n)}n\right)
    =\mu^{(1-\beta)n+o(n)}.
\end{equation*}
By Borel--Cantelli and \eqref{eq:speed-order}, we obtain
$s_{\rm burst}\le s_{\rm cloud}\le1$.
The sustainable speed upper bound follows from the
argument in Section~\ref{sec:sustained}, which applies
to all infinite rays of $\T$.
Finally, the population growth follows from the standard
martingale convergence theorem for Galton--Watson
processes under \eqref{eq:GW2}.
\end{proof}

\section{Convergence in distribution of the maximum}
\label{sec:maximum-distribution}

Write $\gamma=\gamma(m)$. Throughout this section we use the almost-sure limits
\begin{equation}
\lim_{n \rightarrow \infty}\frac{\log_\mu A_n}{n}=\frac{\ainv(M_n)}n=\gamma
\quad\text{almost surely on }\cS.
\label{eq:extreme-growth-input}
\end{equation}
Recall that $\pi_h=\P(Z_h>0)$, and put
$\pi_\infty:=\lim_{h\to\infty}\pi_h>0$.
Since $\pi_\infty\le\pi_h\le1$, the first limit in
\eqref{eq:extreme-growth-input} gives
\begin{equation}
\lim_n\frac{\log_\mu B_n^{\rm eff}}n=\gamma
\qquad\text{almost surely on }\cS,
\label{eq:extreme-effective-bounds}
\end{equation}
where $B_n^{\rm eff}$ is defined in
\eqref{eq:extreme-effective-count}.

For an edge $e=(v,w)$, say that $e$ is
\textit{eligible} at generation $n$ if $v$ is alive
and $w$ has a descendant at generation $n$ in the
unrestricted Galton--Watson tree. Let $\mathcal E_n$
denote the set of eligible edges and define
\begin{equation*}
K_n:=\max_{e\in\mathcal E_n}X_e,
\end{equation*}
with $K_n=0$ when $\mathcal E_n$ is empty.
We first determine the limit of $K_n$, and then
compare it with $M_n$. To begin, we fix
$$
0<c<\frac{\gamma}{\gamma+2},
\qquad k=\lfloor cn\rfloor,
\qquad s:=n-k.
$$
For each $u\in\mathcal A_s$, define
$$
D_{u,k}
:=
\#\bigl\{
(v,w)\in\mathcal E_n:
u\preceq v,\ s<|w|\le n
\bigr\}.
$$
Thus $D_{u,k}$ counts the eligible edges in the
descendant tree of $u$ entering generations
$s+1,\ldots,n$. Put
\begin{equation*}
\Lambda_n
:=
\sum_{u\in\mathcal A_s}
\E[D_{u,k}\mid\cF_s],
\qquad
W_n
:=
\mu\sum_{j=s}^{n-1}\pi_{n-j-1}A_j.
\end{equation*}
Applying branching property and tower property at all generations between $s$ and $n$, we have
$
\E[W_n\mid\cF_s]=\Lambda_n.
$

\begin{lemma}\label{lem:intensity}
\begin{equation}
\lim_{n \rightarrow \infty}\frac{\Lambda_n}{B_n^{\rm eff}}=1,\qquad \text{almost surely on $\cS$}.
\label{eq:extreme-Lambda-Beff}
\end{equation}
\end{lemma}

\begin{proof}
Conditionally on $\cF_s$, the descendant families below the
vertices of $\mathcal A_s$ are independent. The second-moment
bound for the Galton--Watson process gives
\begin{equation}
\operatorname{Var}(W_n\mid\cF_s)
\le C A_s\mu^{2k},
\qquad
\Lambda_n\ge\mu\pi_\infty A_s.
\label{eq:extreme-window-variance}
\end{equation}
Choose $\eta>0$ sufficiently small that
$2c<(\gamma-\eta)(1-c)$ and $(\gamma+\eta)(1-c)<\gamma-\eta.$ 

The first limit in \eqref{eq:extreme-growth-input} implies that
$A_s\ge\mu^{(\gamma-\eta)s}$ eventually on $\cS$.
By Chebyshev and \eqref{eq:extreme-window-variance}, we have
\begin{align*}
&\P\left(\left|W_n-\Lambda_n\right|>
\varepsilon\Lambda_n,\ A_s\ge\mu^{(\gamma-\eta)s}\right)\le
C_\varepsilon\mu^{2k-(\gamma-\eta)s}.
\end{align*}
The right-hand side is summable. Borel--Cantelli therefore yields
\begin{equation}
\lim_{n\rightarrow \infty} \frac{W_n}{\Lambda_n}=1
\qquad\text{almost surely on }\cS.
\label{eq:extreme-W-Lambda}
\end{equation}
It remains to consider the generations before $s$.
Since $\pi_h\le1$ and $\pi_h\ge\pi_\infty$,
\begin{equation*}
0\le B_n^{\rm eff}-W_n
\le\mu\sum_{j<s}A_j,
\qquad
B_n^{\rm eff}\ge\mu\pi_\infty A_{n-1}.
\end{equation*}
By \eqref{eq:extreme-growth-input} and the choice of $\eta$,
$\sum_{j<s}A_j=o(A_{n-1})$ almost surely on $\cS$.
Together with \eqref{eq:extreme-W-Lambda}, this proves the lemma.
\end{proof}

\subsection{The largest displacement}

Write $\Fbar_>(x):=\P(X>x)$. \eqref{eq:tail-regularity} implies
$\Fbar_>(x)\sim \Fbar(x).$
Let $V$ be uniform on $(0,1)$, independent of $X$, and set
$$
U=\Fbar_>(X)+V(\Fbar(X)-\Fbar_>(X)).
$$
Note that $U$ is uniform on $(0,1)$ and
$\Fbar_>(X)\le U\le\Fbar(X)$. Using the above limit, we get
\begin{equation}\label{eq:extreme-tail-rank}
\P(\Fbar(X)\le t)\sim t,\qquad t\downarrow0.
\end{equation}

\begin{lemma}\label{lem:extreme-edge-law}
For every $y>0$, almost surely on $\cS$,
\begin{equation}\label{eq:extreme-conditional-rank}
\lim_{n \rightarrow \infty}\P(\Lambda_n\Fbar(K_n)>y\mid\cF_s)= e^{-y}.
\end{equation}
In particular, under $\P(\,.\mid\cS)$,
\begin{equation}\label{eq:extreme-exp-limit}
\Lambda_n\Fbar(K_n)\xrightarrow{d}E,
\qquad \P(E>y)=e^{-y},\quad y>0.
\end{equation}
\end{lemma}

\begin{proof}
Fix $y>0$ and on $\{\Lambda_n>0\}$ put
$
p_n=\P\left(\Fbar(X)\le y/\Lambda_n\,\middle|\,\cF_s\right),
$
where $X$ is independent of $\cF_s$.
By \eqref{eq:extreme-effective-bounds} and
\eqref{eq:extreme-Lambda-Beff}, $\Lambda_n\to\infty$ on $\cS$.
Thus by \eqref{eq:extreme-tail-rank},
\begin{equation}\label{eq:pn}
\lim_{n \rightarrow \infty}\Lambda_np_n= y
\qquad\text{almost surely on }\cS.
\end{equation}

For $u\in\mathcal A_s$, let $Y_u$ be the number of eligible
edges below $u$ with $\Lambda_n\Fbar(X_e)\le y$, and write
$q_u=\P(Y_u>0\mid\cF_s)$. Since the eligibility of an edge
does not depend on the displacement on that edge,
$
\E[Y_u\mid\cF_s]=p_n\E[D_{u,k}\mid\cF_s].
$
Counting pairs of edges in the unrestricted descendant tree of $u$,
by \eqref{eq:GW2} we have
$$
\E[Y_u(Y_u-1)\mid\cF_s]\le C\mu^{2k}p_n^2.
$$
Using Bonferroni inequalities and $\E[D_{u,k}\mid\cF_s]\le C\mu^k$, and then summing over $u \in \mathcal A_s$ 
we obtain
\begin{align*}
0\le\Lambda_np_n-\sum_{u\in\mathcal A_s}q_u
\le CA_s\mu^{2k}p_n^2 \;\text{ and }\;
\sum_{u\in\mathcal A_s}q_u^2
\le CA_s\mu^{2k}p_n^2.
\end{align*}
We further obtain
$$
A_s\mu^{2k}p_n^2\le C\frac{\mu^{2k}}{A_s}(\Lambda_np_n)^2
=o(1)
\qquad\text{on }\cS,
$$
using \eqref{eq:extreme-window-variance}, \eqref{eq:extreme-growth-input}, the choice of $c$
and \eqref{eq:pn}. Hence almost surely on $\cS$, we have
$\sum_{u\in\mathcal A_s}q_u\rightarrow y$ and $\sum_{u\in\mathcal A_s}q_u^2\rightarrow0$. Therefore, branching property at $s$ yields
\begin{equation}\label{eq:late-edges}
\lim_{n \rightarrow \infty}\P(Y_u=0\text{ for all }u\in\mathcal A_s\mid\cF_s)
=\lim_{n \rightarrow \infty}\prod_{u\in\mathcal A_s}(1-q_u)= e^{-y}
\qquad\text{on }\cS.
\end{equation}

We now consider the displacements up to generation $s$.
Let $J_s^*$ be the largest displacement on an edge entering a generation
at most $s$ whose parent is alive. Note that on $\cS$,
$
J_s^*\le\max\left\{b_s,\max_{1\le j\le s}M_j\right\}.
$
Choose $\delta>0$ such that
$$
\gamma-\delta>1/m,\qquad
(\gamma+\delta)(1-c)<\gamma-\delta,
$$
and put $h_n=a((\gamma-\delta)n)$. By the second limit in
\eqref{eq:extreme-growth-input}, $J_s^*<h_n$ for all sufficiently
large $n$ on $\cS$. Further, by \eqref{eq:extreme-effective-bounds},
\eqref{eq:extreme-Lambda-Beff} and \eqref{C2}, on $\cS$,
$$
\Lambda_n\Fbar(h_n)
=\mu^{\gamma n+o(n)}\mu^{-(\gamma-\delta)n+o(n)}
\longrightarrow\infty.
$$
On the $\cF_s$-measurable event
$\{J_s^*<h_n,\ \Lambda_n\Fbar(h_n)>y\}$,
we thus have $\Lambda_n\Fbar(X_e)>y$ for every eligible edge $e$
with $|e|\le s$. On this event,
$$
\{\Lambda_n\Fbar(K_n)>y\}
=\{Y_u=0\text{ for all }u\in\mathcal A_s\}.
$$
Since this event of occurs eventually on $\cS$,
\eqref{eq:late-edges} provides \eqref{eq:extreme-conditional-rank}. Finally, put $\Omega_n=\{\Lambda_n\Fbar(K_n)>y\}$. We have
\begin{align*}
\left|\P(\Omega_n\cap\cS)
-\E[\1_{\cS}\P(\Omega_n\mid\cF_s)]\right|
\le2\E\left|\1_{\cS}-\P(\cS\mid\cF_s)\right|
\longrightarrow0
\end{align*}
by martingale convergence. Using \eqref{eq:extreme-conditional-rank}
and bounded convergence, we get $\P(\Omega_n\mid\cS)\to e^{-y}$.
\end{proof}

\subsection{Comparison with the alive maximum}

\begin{lemma}\label{lem:extreme-MK-comparison}
Under $\P(\,\cdot\mid\cS)$,
$M_n/K_n\xrightarrow{\P}1.$
\end{lemma}

\begin{proof}
Fix $\eta_0\in(1/2,\gamma)$. The probability that a vertex
at generation $n$ has two ancestral displacements exceeding
$a(\eta_0n)$ is at most
\begin{equation*}
\mu^n\binom n2\Fbar(a(\eta_0n))^2
= n^{O(1)}\mu^{(1-2\eta_0)n+o(n)}.
\end{equation*}
Since $\eta_0>1/2$, by
Borel--Cantelli, every sufficiently long ancestral path contains
at most one displacement exceeding $a(\eta_0n)$.
Each edge on a path ending at an alive vertex of generation $n$
belongs to $\mathcal E_n$. It follows that
\begin{equation*}
M_n\le K_n+n a(\eta_0n)
\qquad\text{eventually on }\cS.
\end{equation*}
Choose $\eta_0<\eta_2<\eta_1<\gamma$. Applying Lemma~\ref{lem:scale-separation}, we have
$a(\eta_1n)>n\,a(\eta_2n)>n^2a(\eta_0n)$ for all sufficiently large $n$. Moreover, by
\eqref{eq:extreme-growth-input}, $M_n\ge a(\eta_1n)$ eventually on $\cS$. Hence
\begin{equation*}
\frac{n a(\eta_0n)}{M_n}\le\frac1n\longrightarrow0
\qquad\text{on }\cS.
\end{equation*}
This proves
$\limsup_{n\to\infty}M_n/K_n\le1$ almost surely on $\cS$.

For the reverse inequality, fix $\varepsilon\in(0,1)$
and $R>0$. Take $h_n=a((\gamma-\delta)n)$ as in the preceding proof and put
$$
H_n:=\left\{
\Lambda_n\Fbar(h_n)>R,\quad
(1-\varepsilon)h_n\ge b_n
\right\}.
$$
The events $H_n$ are $\cF_s$-measurable and occur eventually almost surely on $\cS$.

For each eligible edge $e=(v,w)$ entering a
generation after $s$, choose a path from $w$
to generation $n$ in the unrestricted genealogical
tree, using only the genealogy. Call $e$ bad if
$
\Lambda_n\Fbar(X_e)\le R
$
and some displacement along this path is smaller
than $-\varepsilon h_n/n$.
Put
$
r_n=\P\left(X<-\varepsilon h_n/n\right).
$
Since the chosen path contains at most $n$ edges,
independence of the displacements gives
$$
\P(\text{some bad edge},H_n)
\le Rn r_n.
$$
Indeed, conditionally on $\cF_s$, the expected
number of eligible edges after $s$ is $\Lambda_n$,
while
$$
\P\left(\Lambda_n\Fbar(X)\le R
\,\middle|\,\cF_s\right)
\le\frac{R}{\Lambda_n}.
$$
The probability that a chosen path contains
a displacement smaller than
$-\varepsilon h_n/n$ is at most $nr_n$.
As $h_n\ge b_{n-1}$ eventually,
Assumption~\ref{ass:Lstar} implies that almost surely no bad
edge occurs on $H_n$ for all sufficiently large $n$.

On $H_n\cap\{\Lambda_n\Fbar(K_n)\le R\}$, the edge attaining
$K_n$ enters after $s$ and satisfies $K_n\ge h_n$.
Its chosen descendant path has all subsequent increments at least
$-\varepsilon h_n/n$ and therefore stays above
\begin{equation*}
K_n-\varepsilon h_n
\ge(1-\varepsilon)K_n
\ge b_n.
\end{equation*}
In particular, this path is alive at generation $n$, and
$M_n\ge(1-\varepsilon)K_n$.
By \eqref{eq:extreme-exp-limit},
\begin{equation*}
\limsup_{n\to\infty}
\P(M_n<(1-\varepsilon)K_n\mid\cS)
\le e^{-R}.
\end{equation*}
Letting $R\to\infty$ proves the reverse inequality in probability.
\end{proof}

\begin{proof}[Proof of Theorem~\ref{thm:general-extreme-limit}]
By Lemmas~\ref{lem:intensity} and
\ref{lem:extreme-edge-law},
\begin{equation*}
B_n^{\rm eff}\Fbar(K_n)\xrightarrow{d}E
\end{equation*}
under $\P(\,\cdot\mid\cS)$.
By Lemma~\ref{lem:extreme-MK-comparison}, $M_n/K_n\to1$
in probability. Since $K_n\to\infty$ on $\cS$,
condition \eqref{eq:tail-regularity} gives
\begin{equation*}
\frac{\Fbar(M_n)}{\Fbar(K_n)}\xrightarrow{\P}1.
\end{equation*}
It follows that
$B_n^{\rm eff}\Fbar(M_n)\xrightarrow{d}E.$ Since $\Psi(x)=-\log\Fbar(x)$, taking negative logarithms gives
\begin{equation*}
\Psi(M_n)-\log B_n^{\rm eff}\xrightarrow{d}G,
\qquad
\P(G\le x)=\exp\{-e^{-x}\},\quad x\in\R.
\end{equation*}
This completes the proof.
\end{proof}

\section{Discussion and future directions}\label{sec:discussion}
\begin{enumerate}[label=(\roman*),leftmargin=*]

\item \textit{Growth of the alive population.}
In the present work, we prove the almost sure limit of
$\log_\mu A_n/n$ on $\cS$. It would be interesting 
to investigate whether there exists $a_n, b_n$ and a non-degenerate random variable $V$
such that under $\P(\;.\;| \mathcal S),$ we have
$$\frac{A_n-a_n}{b_n} \xrightarrow{d} V.$$
Using Theorem~\ref{thm:general-extreme-limit}, this would 
provide a better understanding of the fluctuations of
the alive maximum $M_n$.

\item \textit{Extremal point process.}
\citet{LalleySellke1987} proved a conditional limit theorem
for the maximum of branching Brownian motion, which
motivated the study of its extremal point process;
see \citet{ArguinBovierKistler2013,AidekonBerestyckiBrunetShi2013}.
More recently, \citet{YangZhu2024} extended these results
to branching Brownian motion with absorption at a linear
barrier. In our setting, the Poisson approximation obtained
by cutting the tree at a suitable generation in the proof
of Theorem~\ref{thm:general-extreme-limit} also provides
a conditional limit theorem for the maximum. This suggests
studying the convergence of the extremal point process
$$
    \sum_{v\in\mathcal A_n}
    \delta_{\Psi(S(v))-\log B_n^{\mathrm{eff}}},
$$
and understanding how the killing mechanism affects
the genealogical clustering of the extremes.

\item \textit{Genealogy of the alive population.}
Our work has mainly been about understanding the jump structures and how they effect the extremal behaviours. There has not been any focus on the genealogy of the alive particles at generation $n$.

In particular, let $U_n$ and $V_n$
be two typical particles in $\mathcal A_n$, and let $U_n\wedge V_n$ denote
their most recent common ancestor. Can we say anything
about the limiting behaviour of $|U_n\wedge V_n|$
on $\cS$? In particular, does $|U_n\wedge V_n|/n$
converge in distribution ?

\item \textit{Extension of the special-chain argument.}
An extremely interesting observation is that the role of \eqref{eq:BJ}
in the special-chain
argument in Section~\ref{sec:special-chains} is to obtain the bound
$
\tau_{i+1}\le H(e_i)+L,
$
which is used to get the estimate in Lemma~\ref{lem:special-one-step}. What if a weaker bound of the form
$$
\tau_{i+1}\le H(e_i)+o(\tau_{i+1}),
$$
can still be used to obtain a version of
Theorem~\ref{thm:bounded-deficit}.

For instance, consider a displacement law with right tail
$$
\Fbar(x)=\exp\{-(\log x)^\alpha\},
\qquad \alpha>1.
$$
In this case, $a(t)=\exp\{(t\log\mu)^{1/\alpha}\}$
satisfies \eqref{C1}--\eqref{C3} and \eqref{eq:u-l-bound},
but $b_{n-L}/b_n\to1$ for every fixed $L$, and
\eqref{eq:BJ} need not hold. Therefore, under such generalization, one might be able to cover a broader class of distributions and Theorem~\ref{thm:main}
might still hold.

\end{enumerate}
\appendix
\section{Appendix}
\begin{lemma}
\label{lem:scale-separation}
Assume {\rm\eqref{C1}}--{\rm\eqref{C3}}. For every $0<d<c$,
\begin{equation}
\frac{a(ct)}{t\,a(dt)}>1
\qquad\text{for all sufficiently large }t.
\label{eq:scale-separation}
\end{equation}
Consequently, for every $\beta>0$, $B<m\beta$
and $K<\infty$, there exists $t_0<\infty$ such that
\begin{equation}
a(\beta t)-K\lceil Bt\rceil
\ge b_{\lceil Bt\rceil},
\qquad t\ge t_0.
\label{eq:barrier-buffer}
\end{equation}
\end{lemma}

\begin{proof}
Suppose that \eqref{eq:scale-separation} fails.
Then there exist arbitrarily large $t$ such that $a(ct)/t\le a(dt).$ Since $\Fbar$ is non-increasing, this implies
\begin{equation*}
\Fbar\left(\frac{a(ct)}{t}\right)
\ge \Fbar(a(dt)).
\end{equation*}
By Assumption {\rm\eqref{C3}}, applied with $r=1/c$,
\begin{equation*}
\log_\mu\Fbar\left(\frac{a(ct)}{t}\right)
=-ct+o(t).
\end{equation*}
On the other hand, Assumption {\rm\eqref{C2}} gives
$\log_\mu\Fbar(a(dt))=-dt+o(t).$
Since $c>d$, these two estimates contradict the
preceding inequality for all sufficiently large $t$.
This proves \eqref{eq:scale-separation}.

For the second assertion, fix $\beta>0$, $B<m\beta$, $K<\infty$ and choose $d \in (B/m, \beta)$. 
For all sufficiently large $t$,
$\lceil Bt\rceil/m\le dt$, and hence
$b_{\lceil Bt\rceil}
=a\left(\lceil Bt\rceil/m\right)
\le a(dt).$

By \eqref{eq:scale-separation},
$a(\beta t)>t\,a(dt)$ for all sufficiently large $t$.
Moreover, Assumption {\rm\eqref{C1}} implies that
$a(dt)\to\infty$, and therefore
$(t-1)a(dt)\ge K\lceil Bt\rceil$
for all sufficiently large $t$. Consequently,
\begin{align*}
a(\beta t)>t\,a(dt)=a(dt)+(t-1)a(dt)
\ge b_{\lceil Bt\rceil}+K\lceil Bt\rceil.
\end{align*}
This proves \eqref{eq:barrier-buffer}.
\end{proof}

\begin{proof}[Proof of Lemma~\ref{lem:examples}]
By the choice of $a$, for all sufficiently large $t$,
\begin{equation}
\Fbar(a(t))=\mu^{-t}.
\label{eq:canonical-u-l-bound}
\end{equation}
Thus \eqref{C2} and \eqref{eq:u-l-bound} hold, the latter with
$c_*=C_*=1$. We verify \eqref{C1} and \eqref{C3} for the two classes.

\emph{Case (i).} Define
$U(y):=\inf\{x>0:\Fbar(x)\le y^{-1}\}$. Since
$\Fbar\in RV_{-\alpha}$, asymptotic inversion gives
\begin{equation*}
U\in RV_{1/\alpha},\qquad a(t)=U(\mu^t),\qquad
\log a(t)=\frac{t\log\mu}{\alpha}+o(t).
\end{equation*}
In particular, \eqref{C1} holds. Moreover, for fixed $s\in\mathbb R$,
\begin{equation}
\frac{a(t+s)}{a(t)}\longrightarrow\mu^{s/\alpha}.
\label{eq:RV-additive-shift}
\end{equation}
Fix $r>0$. By Potter's bounds, for every $\varepsilon>0$ and
all sufficiently large $t$,
\begin{equation*}
C_\varepsilon^{-1}t^{\alpha-\varepsilon}
\le\frac{\Fbar(a(t/r)/t)}{\Fbar(a(t/r))}
\le C_\varepsilon t^{\alpha+\varepsilon}.
\end{equation*}
Together with \eqref{eq:canonical-u-l-bound}, this implies
$\log_\mu\Fbar(a(t/r)/t)=-t/r+o(t)$, proving \eqref{C3}.

\emph{Case (ii).} Put $W:=V^{-1}$. Then
\begin{equation*}
\log a(t)=W(t\log\mu),\qquad W\in RV_{1/\beta}.
\end{equation*}
Since $1/\beta>1$, \eqref{C1} follows. For fixed $r>0$, let
$z_t:=W(t\log\mu/r)$. We have $z_t/\log t\to\infty$, so regular
variation of $V$ gives
\begin{align*}
-\log\Fbar(a(t/r)/t)
&=V(z_t-\log t)
 =V(z_t)(1+o(1))
 =\frac{t\log\mu}{r}(1+o(1)).
\end{align*}
This proves \eqref{C3}.

In both cases there is $c>0$ such that $b_n\ge e^{cn}$ for all
sufficiently large $n$. Therefore, for every $c_1>0$, the
left-tail bound \eqref{eq:polynomial-left-tail} yields
\begin{equation*}
\P\left(X<-c_1\frac{b_{n-1}}n\right)
\le C_-c_1^{-\eta}n^\eta b_{n-1}^{-\eta}
\le Ce^{-c'n}
\end{equation*}
for some $c'>0$. Hence Assumption~\ref{ass:Lstar} holds.

It remains to verify \eqref{eq:BJ}. Set
$c_0:=(\sum_{j\ge1}j^{-2})^{-1}$, choose an integer
$r_0>\max\{m,1\}$, and then $\delta>0$ such that
\begin{equation}
\delta<\frac{c_0}{(r_0-1)^2}.
\label{eq:BJ-delta}
\end{equation}
We claim that for a suitable integer $L_{\rm BJ}\ge1$,
$b_{n-L_{\rm BJ}}\le\delta b_n$ for all sufficiently large $n$. In case (i),
\eqref{eq:RV-additive-shift} gives
$b_{n-L}/b_n\to\mu^{-L/(m\alpha)}$, and we choose $L$ large.
In case (ii), \eqref{eq:LW-smooth} implies
\begin{equation*}
W'(s)=\frac1{V'(W(s))}\sim\frac{W(s)}{\beta s}
\longrightarrow\infty.
\end{equation*}
Thus $b_{n-1}/b_n\to0$ by the mean value theorem, so
$L_{\rm BJ}=1$ is sufficient in this case.

Let $Y_i=X_i^+$ and write $Y_{(1)}\ge\cdots\ge Y_{(n)}$
for the decreasing rearrangement. Whenever $\sum_iY_i\ge x$,
there is some $k\le n$ such that
\begin{equation}
Y_{(k)}\ge\frac{c_0x}{k^2};
\label{eq:order-stat}
\end{equation}
otherwise the sum would be smaller than
$c_0x\sum_{k\ge1}k^{-2}=x$. On the event
\begin{equation*}
\sum_{i=1}^nX_i\ge b_n,
\qquad \max_{i\le n}X_i<b_{n-L_{\rm BJ}},
\end{equation*}
we have $\sum_iY_i\ge b_n$ and $Y_{(1)}<\delta b_n$.
By \eqref{eq:BJ-delta}, the index in
\eqref{eq:order-stat} must satisfy
$k\ge r_0$. Consequently,
\begin{align}
&\P\left(\sum_{i=1}^nX_i\ge b_n,
\ \max_{i\le n}X_i<b_{n-L_{\rm BJ}}\right)
\le
\sum_{k=r_0}^n{n\choose k}
\Fbar\left(\frac{c_0b_n}{k^2}\right)^k.
\label{eq:BJ-union}
\end{align}

Choose $\lambda$ with $1/r_0<\lambda<1/m$. We claim that
\begin{equation}
\sup_{1\le k\le n}
\Fbar\left(\frac{c_0b_n}{k^2}\right)
\le\mu^{-\lambda n}
\label{eq:BJ-uniform-tail-bound}
\end{equation}
for all sufficiently large $n$. In case (i), Potter's bounds give
$\Fbar(c_0b_n/k^2)\le Cn^d\Fbar(b_n)
=Cn^d\mu^{-n/m}$ uniformly for $k\le n$, which proves the claim.
In case (ii), $\log b_n/\log n\to\infty$, and hence
\begin{equation*}
\sup_{1\le k\le n}
\left|\frac{\log(c_0b_n/k^2)}{\log b_n}-1\right|
\longrightarrow0.
\end{equation*}
Regular variation of $V$ and $V(\log b_n)=n\log\mu/m$
then give
\begin{equation*}
\sup_{1\le k\le n}
\Fbar\left(\frac{c_0b_n}{k^2}\right)
\le\mu^{-(1/m-o(1))n},
\end{equation*}
which proves \eqref{eq:BJ-uniform-tail-bound}. By \eqref{eq:BJ-union},
${n\choose k}\le(en/k)^k$, and
$q_n:=en\mu^{-\lambda n}\to0$,
\begin{align*}
&\mu^n\P\left(\sum_{i=1}^nX_i\ge b_n,
\ \max_{i\le n}X_i<b_{n-L_{\rm BJ}}\right)
\le\mu^n\sum_{k=r_0}^\infty q_n^k
\le2(en)^{r_0}\mu^{(1-r_0\lambda)n}.
\end{align*}
Since $r_0\lambda>1$, the last expression is summable in $n$. This proves \eqref{eq:BJ}.

Finally, we verify \eqref{eq:tail-regularity}. In case~(i), since $L$ is slowly varying,
\begin{equation*}
\frac{\Fbar((1-\varepsilon)x)}{\Fbar((1+\varepsilon)x)}\longrightarrow\Bigl(\frac{1+\varepsilon}{1-\varepsilon}\Bigr)^{\alpha},\qquad x\to\infty,
\end{equation*}
and the right-hand side tends to $1$ as $\varepsilon\downarrow0$. In case~(ii), \eqref{eq:LW-smooth} gives $V'(y)\sim\beta V(y)/y\to0$, since $\beta<1$. Hence, by the mean value theorem,
\begin{equation*}
\log\frac{\Fbar((1-\varepsilon)x)}{\Fbar((1+\varepsilon)x)}
=V(\log x+\log(1+\varepsilon))-V(\log x+\log(1-\varepsilon))\longrightarrow0
\end{equation*}
for each fixed $\varepsilon\in(0,1)$. This proves \eqref{eq:tail-regularity} in both cases and completes the proof.
\end{proof}
\paragraph{Acknowledgement}
Ramkrishna Jyoti Samanta acknowledges support from the Additional Funding Programme for Mathematical Sciences, delivered by EPSRC (EP/V521917/1)
and the Heilbronn Institute for Mathematical Research. The authors also acknowledge the EPSRC funded workshop "Probability meets biology iii" at the University of Bath, where the authors met and discussed their progress on this project.
\paragraph{Statements and Declarations}
The authors declare no competing interests. There is no use of data in this work, so data availability is not applicable here. The authors themselves are accountable for this article. 

\end{document}